\documentclass{article}

\usepackage[final,nonatbib]{nips_2018}

\usepackage[utf8]{inputenc} % allow utf-8 input}
\usepackage[T1]{fontenc}    % use 8-bit T1 fonts
\usepackage{hyperref}       % hyperlinks
\usepackage{url}            % simple URL typesetting
\usepackage{booktabs}       % professional-quality tables
\usepackage{amsfonts}       % blackboard math symbols
\usepackage{nicefrac}       % compact symbols for 1/2, etc.
\usepackage{microtype}      % microtypography

\usepackage[dvipdfmx]{graphicx} 
\usepackage{bmpsize}

\usepackage{mathtools}
\DeclarePairedDelimiter{\ceil}{\lceil}{\rceil}

\usepackage{amsmath, amsfonts, amssymb, amsthm}
\usepackage{xcolor}
\usepackage{algorithm,algorithmic}
\def\W{\mathcal W}
\def\C{\mathcal C}
\def\X{\mathcal X}
\def\Y{\mathcal Y}
\def\Z{\mathcal Z}
\def\M{\mathcal M_+^1}
\def\R{\mathbb R}
\def\E{\mathbb E}

\def\Blm{\boldsymbol{\lambda}}

\def\Bzeta{\boldsymbol{\zeta}}
\def\Beta{\boldsymbol{\eta}}
\def\BBlm{\boldsymbol{\bar{\lambda}}}

\def\BBzeta{\boldsymbol{\bar{\zeta}}}
\def\BBeta{\boldsymbol{\bar{\eta}}}
\def\p{\mathtt{p}}
\def\blm{\bar{\lambda}}

\def\lm{\lambda}
\def\eig{\lambda_{\max}(W)}

\newcommand{\e}{\varepsilon}
\newcommand{\la}{\langle}
\newcommand{\ra}{\rangle}
\def\vp{\varphi}

\newtheorem{theorem}{Theorem}
\newtheorem{lemma}{Lemma}

\def\cu{}%#1{{\color{red}#1}}  % Cesar
\def\dd{}%#1{{\color{magenta}#1}} %Darina Dvinskikh's text
\usepackage{times} % assumes new font selection scheme installed
\usepackage{amsmath} % assumes amsmath package installed
\usepackage{amssymb}
\usepackage{amsfonts}  % assumes amsmath package installed
\usepackage{latexsym}
\usepackage{epsfig}
\usepackage{graphicx} 
\usepackage{dsfont}
\usepackage{mathtools}
\usepackage{subfigure}
\usepackage{lscape}
\usepackage{color}
\usepackage{enumerate}
\usepackage[percent]{overpic}
\usepackage{hyperref}
\usepackage[export]{adjustbox}
\usepackage{tikz}
\usepackage{enumitem}

\usepackage{pgfplotstable}
\usepackage{pgfplots}
\pgfplotsset{
	table/search path={plot_figures},
}
\usepackage{tikz}
\usepackage{xr}
\usetikzlibrary{external}
\usetikzlibrary{arrows,%
	petri,%
	topaths}%
\usetikzlibrary{graphs,graphs.standard}
\usepackage{tkz-berge}
\usepackage{subfigure}
\pgfdeclarelayer{bg}    % declare background layer
\pgfsetlayers{bg,main}  % set the order of the layers (main is the standard layer)
\usepgfplotslibrary{polar}
\pgfplotsset{compat=1.14}

\graphicspath{{Figures/}}

\title{Decentralize and Randomize: Faster Algorithm for Wasserstein Barycenters}%Fast Stochastic Distributed Convex Optimization for Wasserstein Barycenters over Networks}

\author{
Pavel Dvurechensky, Darina Dvinskikh \\
Weierstrass Institute for Applied Analysis and Stochastics, \\
Institute for Information Transmission Problems RAS\\
\texttt{\{pavel.dvurechensky,darina.dvinskikh\}@wias-berlin.de}
\And
Alexander Gasnikov \\
Moscow Institute of Physics and Technology,\\ Institute for Information Transmission Problems RAS \\
\texttt{gasnikov@yandex.ru} 
\And
C\'{e}sar A. Uribe \\
Massachusetts Institute of Technology \\
\texttt{cauribe@mit.edu}
\And
Angelia Nedi\'{c} \\
Arizona State University, \\ Moscow Institute of Physics and Technology \\
\texttt{angelia.nedich@asu.edu}
}

\begin{document}

\maketitle
\begin{abstract}
We study \cu{the} decentralized distributed computation of discrete approximations for the regularized Wasserstein barycenter of a finite set of continuous probability measures distributedly stored over a network. We assume there is a network of agents/machines/computers, \cu{and} each agent holds a private continuous probability measure and seeks to compute the barycenter of all the measures in the network by getting samples from its local measure and exchanging information with its neighbors. Motivated by this problem, we develop, and analyze, a novel accelerated primal-dual stochastic gradient method for general stochastic convex optimization problems with linear equality constraints. Then, we apply this method to the decentralized distributed optimization setting to obtain a new algorithm for the distributed semi-discrete regularized Wasserstein barycenter problem. \cu{Moreover},  we show explicit non-asymptotic complexity for the \cu{proposed} algorithm. 
%The proposed algorithm can be executed over arbitrary networks that are undirected, connected and static, using the local information only. 
%Moreover, we show explicit non-asymptotic complexity in terms of the problem parameters. 
Finally, we show the effectiveness of our method on the distributed computation of the regularized Wasserstein barycenter of univariate Gaussian and von Mises distributions, as well as some applications to image aggregation.
\end{abstract}

\section{Introduction}
Optimal transport (OT) \cite{Monge1781,kantorovich1942translocation} has become increasingly popular in the machine learning and optimization community. 
Given a basis space (e.g., pixel grid) and a transportation cost function (e.g., squared Euclidean distance), the OT approach defines a distance between two objects (e.g., images), modeled as two probability measures on the basis space, as the minimal cost of transportation of the first measure to the second.
%This approach allows defining a distance between two objects, modeled as two probability measures, as the minimal cost of transportation of the first measure to the second.
Besides images, these probability measures or histograms can model other real-world objects like videos, texts, etc. The optimal transport distance leads to outstanding results in unsupervised learning \cite{arjovsky2017wasserstein,bigot2017geodesic}, semi-supervised learning \cite{solomon2014wasserstein}, clustering \cite{ho17multilevel}, text classification \cite{kusner2015from}, as well as in image retrieval, clustering and classification \cite{rubner2000earth,cuturi2013sinkhorn,sandler2011nonnegative}, statistics \cite{ebert2017construction,panaretos2016amplitude}, economics and finance~\cite{Beiglbock2013}, condensed matter physics~\cite{Buttazzo2012}, and other applications \cite{kolouri2017optimal}. From the computational point of view, the optimal transport distance (or Wasserstein distance) between two histograms of size $n$ requires solving a linear program, which typically requires $O(n^3\log n)$ arithmetic operations. An alternative approach is based on entropic regularization of this linear program and application of either Sinkhorn's algorithm \cite{cuturi2013sinkhorn} or stochastic gradient descent \cite{genevay2016stochastic}, both requiring $O(n^2)$ arithmetic operations, which can be too costly in the large-scale context. 

Given a set of objects, the optimal transport distance naturally defines their mean representative. \cu{For example}, the $2$-Wasserstein barycenter \cite{Agueh2011} \cu{is} an object minimizing the sum of squared $2$-Wasserstein distances to all objects in \cu{a} set. Wasserstein barycenters capture the geometric structure of objects, such as images, better than the barycenter with respect to the Euclidean or other distances \cite{Cuturi2014}. If the objects in the set are randomly sampled from some distribution, theoretical results such as central limit theorem \cite{barrio1999central} or confidence set construction \cite{ebert2017construction} have been proposed, providing the basis for the practical use of Wasserstein barycenter. However, %the computational aspects present some difficulties: 
calculating the Wasserstein barycenter of $m$ measures includes repeated \cu{computation} of $m$ Wasserstein distances. The entropic regularization approach was extended for this case in \cite{benamou2015iterative}, with the proposed algorithm having \cu{a $O(mn^2)$ complexity}, which \cu{can be} very large if $m$ and $n$ are large. Moreover, in the large-scale setup, storage and processing of transportation plans, required to calculate Wasserstein distances, can \cu{be intractable for local computation.}
%lead to an overflow of computer's memory.
%no single computer can store transport plans, required for calculating the Wasserstein distances, in RAM. 
On the other hand, recent studies \cite{Nedic2017achieving, Scaman2017Optimal, Optimal_Varying,uribe2018dual,Nedic2017cc} on accelerated distributed convex optimization algorithms demonstrated their efficiency for convex optimization problems over arbitrary networks with inherently distributed data, \cu{i.e.,} the data is produced by a distributed network of sensors \cite{olfati-saber2006belief,Nedic2017bb,Nedic2017aa} or the transmission of information is limited by communication or privacy constraints, i.e., only limited amount of information can be shared across the network. %meaning that the agents share neither their local functions nor their gradient information explicitly.

Motivated by the limited communication issue and the computational complexity of the Wasserstein barycenter problem for large amounts of data stored in a network of computers, we use the entropy regularization of the Wasserstein distance and propose a decentralized algorithm to calculate an approximation to the  Wasserstein barycenter of a set of probability measures. We solve the problem in a distributed manner on a connected and undirected network of agents oblivious to the network topology. Each agent locally holds a possibly continuous probability distribution, can sample from it, and seeks to cooperatively compute the barycenter of all probability measures exchanging the information with its neighbors. We consider the semi-discrete case, which means that we fix the discrete support for the barycenter and calculate a discrete approximation for the barycenter.

\subsection{Related work}
Unlike \cite{uribe2018distributed}, we propose a decentralized distributed algorithm for the computation of the regularized Wasserstein barycenter of a set of \textit{continuous} measures. Working with continuous distributions requires the application of stochastic procedures like stochastic gradient method as in \cite{genevay2016stochastic}, where it is applied for regularized Wasserstein distance, but not for Wasserstein barycenter. This idea was extended to the case of non-regularized barycenter in \cite{Staib2017,Claici2018StochasticWB}, where parallel algorithms were developed. The important difference between the parallel and the decentralized setting is that, in the former, the topology of the computational network is fixed to be a star topology and it is known in advance by all the machines, forming a master/slave architecture. We seek to further scale up the barycenter computation to a huge number of input measures using \textit{arbitrary} network topologies. Moreover, unlike \cite{Staib2017}, we use entropic regularization to take advantage of the problem smoothness and obtain faster rates of convergence for the optimization procedure. Unlike \cite{Claici2018StochasticWB}, we fix the support of the barycenter, which leads to a convex optimization problem and allows us to prove complexity bounds for our algorithm. 

The well-developed approach based on Sinkhorn's algorithm \cite{cuturi2013sinkhorn,benamou2015iterative,cuturi2016smoothed} naturally leads to parallel algorithms. Nevertheless, its application to continuous distributions requires discretization of these distributions, leading to computational intractability when one desires good accuracy and, hence, has to use fine discretization with large $n$, which leads to the necessity of solving an optimization problem of large dimension. Thus, this approach is not directly applicable in our setting of continuous distributions, and it is not clear whether it is applicable in the decentralized distributed setting with arbitrary networks.

Recently, an alternative accelerated-gradient-based approach was shown to give better results than the Sinkhorn's algorithm for Wasserstein distance \cite{dvurechensky2018computational,dvurechensky2017adaptive}. Moreover, accelerated gradient methods have natural extensions for the decentralized distributed setting \cite{Scaman2017Optimal,uribe2017optimal,lan2017communication}.
Nevertheless, existing distributed optimization algorithms can not be applied to the barycenter problem in our setting of continuous distributions as these algorithms are either designed for deterministic problems or for stochastic primal problem, whereas in our case the \emph{dual} problem is a stochastic problem.
Table \ref{T:PapCompare} summarizes the existing literature on Wasserstein barycenter calculation and shows our contribution.

%Problem of calculation of Wasserstein distance between discrete measures is turned to linear program~\cite{Anderes2016}, in case of continuous measures or discrete and continuous the Wasserstein distance between them requires stochastic approaches \cite{Genevay2016}. Due to the high complexity of computation Wasserstein distance also known as Earth Mover's Distance a lot of works \cite{Cuturi2014, Cuturi2016, Genevay2016} smooth Wasserstein distance by adding an entropic regularizer to reduce the computation complexity.\\  The authors of \cite{Staib2017} use a parallelizable structure of dual problem of computing the Wasserstein barycenter. They also work with arbitrary measure and discretize only their barycenter. However, parallel technique concedes decentralized technique for large-scale problems because it requires a master center that coordinates the actions of the parallel machines.

\vspace{-0.2cm}
\begin{table}[h]
\caption{Summary of our contribution.}
 \vspace{-0.5cm}
\label{T:PapCompare}
\vskip 0.15in
\begin{center}
\begin{small}
\begin{sc}
\begin{tabular}{lcccr}
\toprule
Paper & Decentralized & Continuous & Barycenter \\
\midrule
%Cuturi \cite{cuturi2013sinkhorn}, Benamou et.al. \cite{benamou2015iterative}, Cuturi, Peyr\'e \cite{cuturi2016smoothed}  & $\times$ & $\times$ & $\surd$ \\
\cite{cuturi2013sinkhorn,benamou2015iterative,cuturi2016smoothed}  & $\times$ & $\times$ & $\surd$ \\
% & $\times$ & $\times$ & $\surd$ \\
\cite{genevay2016stochastic} & $\times$ & $\surd$ & $\times$ \\
\cite{Staib2017,Claici2018StochasticWB}  & $\times$ & $\surd$ & $\surd$ \\
%Distributed Computation of Wasserstein Barycenters over
%Networks \\ (not sure that we should mention it) & $\surd$ & $\times$ & $\surd$ \\
This paper (Alg. \ref{alg:main}) & $\surd$ & $\surd$ & $\surd$  \\
\bottomrule
\end{tabular}
\end{sc}
\end{small}
\end{center}
\vskip -0.1in
\end{table}
\vspace{-0.2cm}

%\footnotetext{They consider a parallel algorithm which is a very particular case of distributed optimization.}


% acknowledgments go at the end of the paper. Do not include
% acknowledgments in the anonymized submission, only in the final paper.

\subsection{Contributions}
%Initially, we formulate the Wasserstein barycenter problem as a convex distribution optimization problem with a set of consensus constraints. These constraints are induced by the topology of the communication network where agents interact and are represented as a set of linear equality constraint. Then, we construct the dual problem and show that it is a smooth stochastic optimization problem, which can be solved in a distributed manner with agents communicating only with their neighbors without knowing the whole topology of the network. We solve the dual problem and reconstruct the solution of the primal problem, i.e., the barycenter, by proposing a novel algorithm for distributed stochastic optimization, the Accelerated Primal-Dual Stochastic Gradient Descent method (APDSGD). We summarize our contributions as follows:

\begin{itemize}[leftmargin=*]
\item We propose a novel algorithm for general stochastic optimization problems with linear constraints, namely the Accelerated Primal-Dual Stochastic Gradient Method (APDSGM).
\item We propose a distributed algorithm for the computation of a discrete approximation for regularized Wasserstein barycenter\cu{s} of a set of continuous distributions stored distributedly over a network (connected and undirected) with unknown arbitrary topology.
\item We provide iteration and arithmetic operations complexity for the proposed algorithms in terms of the problem parameters.%, the network topology, the size of the discretization of the barycenter and the desired quality of the solution.
\item  We demonstrate the effectiveness of our algorithm on the distributed computation of the regularized Wasserstein barycenter of a set of Gaussian distributions and a set of von Mises distributions for various network topologies and network sizes. Moreover, we show some initial results on the problem of image aggregation for two datasets, namely, a subset of the MNIST digit dataset~\cite{LeCun1998} and subset of the IXI Magnetic Resonance dataset~\cite{ixidata}. %\item Complexity bounds for the problem of 
\end{itemize}

%1. Our main contribution is presenting regularized Wasserstein barycenter for semi-discrete optimal transport in the form which can be executed in a distributed manner for further distributed calculation.  
%2.We introduce a stochastic algorithm executed in a distributed manner estimating the discrete approximation of smoothed Wasserstein barycenter of a set of arbitrary measures using mini-batch for better unbiased estimation of the gradient in gradient descent step. \\
%2. We estimate the size of mini-batch for calculating the barycenter with the desired accuracy.\\

\subsection{Paper organization}
This paper is organized as follows. In Section~\ref{sec:problem}, we present the regularized Wasserstein barycenter problem for the semi-discrete case and its distributed computation over networks. In Section~\ref{sec:stoch_prob}, we introduce a new algorithm for general stochastic optimization problems with linear constraints and analyze its convergence rate. Section \ref{sec:GenAlgApplication} extends this algorithm and introduces our method for the distributed computation of regularized Wasserstein barycenter. Section \ref{sec:experiments} shows the experimental results for the proposed algorithm. The appendix contains the proofs of stated lemmas and theorems, as well as additional results of numerical experiments.

\textbf{Notation.} We define $\M(\X)$ -- the set of positive Radon probability measures on a metric space $\X$, and $S_1(n)  = \{ a \in \mathbb{R}_+^n  \mid \sum_{l=1}^n a_l =1 \}$ the probability simplex. We use $\C(\X)$ as the space of continuous functions on $\X$. We denote by $\delta(x)$ the Dirac measure at point $x$. We refer to $\lambda_{\max}(W)$ as the maximum eigenvalue of matrix W. We also use bold symbols for stacked vectors $\p = [p_1^T,\cdots,p_m^T]^T \in \mathbb{R}^{mn}$, where $p_1,...,p_m\in \R^n$. In this case $[\p]_i=p_i$ -- the $i$-th block of $\p$. For a vector $\lambda \in \R^n$, we denote by $[\lambda]_l$ its $l$-th component. We refer to the Euclidean norm of a vector $\|p\|_2:=\sum_{l=1}^n([p]_l)^2$ as $2$-norm. 

%Mathematical expectation is denoted by $\mathbb{E}$ . 
%We define the Kronecker product $\otimes$.

\section{The Distributed Wasserstein Barycenter Problem}\label{sec:problem}

In this section, we present the problem of decentralized distributed computation of regularized Wasserstein \cu{barycenters} for a family of possibly continuous probability measures distributed over a network. 
First, we provide the necessary background for %regularized Wasserstein distance and barycenter. 
entropic regularization of optimal transport and \cu{the} Wasserstein barycenter problem.
Then, we give the details of the distributed formulation of the optimization problem defining \cu{the} Wasserstein barycenter, which is a minimization problem with linear equality constraint. To deal with this constraints, we make a transition to the dual problem, which, as we show, due to the presence of continuous distributions, is a smooth stochastic optimization problem. 
%Next, we will state this problem as a couple of primal and dual stochastic optimization problems. 

%This sharing constraint can be written as a system of linear equations and the minimization  problem for finding the barycenter becomes a linearly constrained optimization problem. A natural way to solve it, is to form the dual problem. Then each agent is possessing a dual variable and a part of the dual objective which turns out to be smooth. In our setting, for each agent, its part of the dual objective can be interpreted as an expectation of some function with respect to the measure corresponding to this agent. To sum up, the dual problem is  a smooth stochastic optimization problem. 
%\pd{TODO: clarify the distinction between measure and probability distribution. Is it ok to say that we sample from a measure?}
%\dd{REMARK: Any positive distribution defines a positive Radon measure. \\
%It is ok to say "we sample from a measure"(see i.e. M.Cuturi, M.Staib (Parallel Streaming...))}. 

\subsection{Regularized semi-discrete formulation of the optimal transport problem}
%Assume that we are given a positive Radon probability measure $\mu$  on the space $\mathcal{Y}$ with density $q(y)$  with respect to some reference measure $dy$ and a discrete probability measure $\nu=\sum_{i=1}^np_i\delta(z_i)$ on the space $\mathcal{Z}$ with reference measure $dz$. We also assume that the reference measure on $\mathcal{Y} \times \mathcal{Z}$ is given by $dydz$
%We consider Wasserstein distance and its entropic regularization in the semi-discrete setting. 
%\begin{equation}
%\W_{\gamma}(\mu,\nu):=\min\left\{ \int_{\mathcal{Y}} \sum_{i=1}^n (c(y,z_i)   + \gamma  \ln ( \pi_i(y) ) ) \pi_i(y)dy: \sum_{i=1}^n \pi_i(y) = q(y), y \in \mathcal{Y}; \int_{\mathcal{Y}} \pi_i(y)dy = p_i, i=1,...,n  \right\}
%\end{equation}
We consider entropic regularization for the optimal transport problem and the corresponding regularized Wasserstein distance and barycenter~\cite{cuturi2013sinkhorn}. Let  $\mu \in \mathcal{M}_{+}^1(\mathcal{Y})$  with density $q(y)$  on a metric space $\mathcal{Y}$, and a discrete probability measure $\nu=\sum_{i=1}^n[p]_i\delta(z_i)$ with weights given by vector $p \in S_1(n)$ and finite support given by points $z_1, \dots, z_n \in \mathcal{Z}$  from a metric space $\mathcal{Z}$. Denote by $c_i(y) = c(z_i,y)$ a cost function for transportation of a unit of mass from point $z_i \in \mathcal{Z}$  to point $y \in \mathcal{Y}$. 
%We consider a Kantorovich formulation of the optimal transport problem with Kullback-Leibler regularization. 
Then we define regularized Wasserstein distance in semi-discrete setting between continuous measure $\mu$ and discrete measure $\nu$ as follows\footnote{Formally, the $\rho$-Wasserstein distance for $\rho \geq 1$ is $\left({\W_0(\mu,\nu)}\right)^{\frac{1}{\rho}}$ if $\Y=\Z$ and $c_i(y) = d^{\rho}(z_i,y)$, $d$ being a distance on $\Y$. 
For simplicity, we refer to \eqref{WassDis} as regularized Wasserstein distance in a general situation since our algorithm does not rely on any specific choice of cost $c_i(y)$.}  
\begin{align}\label{WassDis}
\W_{\gamma}(\mu,\nu)=\min_{\pi \in \Pi(\mu,\nu)}\left\{\sum_{i=1}^n\int_{\Y}  c_i(y)\pi_i(y)dy   + \gamma \sum_{i=1}\limits^n\int_{\Y}  \pi_i(y)\log\left(\frac{\pi_i(y)}{\xi}\right)dy\right\},
\end{align}
where $\xi$ is the uniform distribution on $\Y \times \mathcal{Z}$, 
%$d\xi(i,z) = \frac{dz}{n\int_\Y dy}$ $(z\in \Y, ~i=1,\dots,n)$  and 
%{\sloppy $KL(\pi|\xi)=$}, 
and the set of admissible coupling measures $\pi$ is defined as follows
\begin{align*}
\Pi(\mu,\nu) = \left\{\pi \in\M(\Y) \times  S_1(n): \sum_{i=1}^n \pi_i(y) = q(y), y \in \mathcal{Y}, \int_{\mathcal{Y}} \pi_i(y)dy = p_i, \forall~i=1,\dots,n  \right\}.   
\end{align*}
We emphasize that, unlike \cite{genevay2016stochastic}, we regularize the problem by the Kullback-Leibler divergence from the uniform distribution $\xi$, which allows us to find explicitly the Fenchel conjugate for $\W_{\gamma}(\mu,\nu)$, see Lemma \ref{Lm:dual_obj_properties} below.

%use another measure $\xi$ in unlike in, where $\xi$ is a product of measure $\mu\otimes \nu$. The goal of using univariate measure $\xi$ is further separating dual variables in the dual problem as it was done in \cite{Uribe2018}, where the decentralized technique was applied to discrete optimal transport. 

%\pd{I think it is interesting to cover also the non-regularized setting. I'm not sure that for this case the regularization should look like this. Note that \url{https://arxiv.org/pdf/1605.08527.pdf} use another regularization, which, I think leads to another dual problem.}\\

For a set of measures $\mu_i \in \M(\mathcal{Z})$, $i=1, \dots, m$, we fix the support  $z_1, \dots, z_n \in \mathcal{Z}$ of their regularized Wasserstein barycenter $\nu$ and wish to find it in the form $\nu=\sum_{i=1}^n[p]_i\delta(z_i)$, where $p \in S_n(1)$. %The motivation behind considering fixed support is an intention to get the convex problem and better theoretical estimations. 
Then the regularized Wasserstein barycenter in the semi-discrete setting is defined as the solution to the following convex optimization problem\footnote{For simplicity, we assume equal weights for each $\W_{\gamma,\mu_i}(p)$ and do not normalize the sum dividing by $m$. Our results can be directly generalized to the case of non-negative weights summing up to 1.}
%Here we suppose the distance for each measure to be present in the sum with the same weights, for simplicity, equal to 1 but not $1/m$, as it is often considered in the distributed optimization.  }
%For a set of positive Radon probability measures $(\mu_1, \dots, \mu_m)$ the regularized Wasserstein barycenter in the semi-discrete setting is defined as the solution to the following convex optimization problem\footnote{For simplicity, we consider the distance for each measure to be present in the sum with the same weight equal to 1. Our results can be directly generalized to the case of different weights and equal weights e.g. $1/m$ as it often encounters in the literature.}
\begin{align}\label{w_barycenter}
\min_{p \in S_1(n)} \sum\limits_{i=1}^{m}  \W_{\gamma,\mu_i}(p),
\end{align}
where we 
%fixed the support $z_1, \dots, z_n \in \mathcal{Z}$ of the barycenter $\nu$ and characterize it by the vector $p \in S_n(1)$, i.e., $\nu=\sum_{i=1}^np_i\delta(z_i)$ and 
used notation $\W_{\gamma,\mu}(p) := \W_{\gamma}(\mu,\nu)$ for fixed probability measure $\mu$.

%\subsection{Wasserstein distance and barycenter in semi-discrete formulation}
%The Wasserstein distance corresponds to the dual formulation of semi-discrete OT is given by
%\begin{align}\label{dual_OT}
%\mathcal{W}_\gamma(\mu,\nu) = \max_{(g,\lambda)} \left\{ \int_\mathcal{Y}g(y)q(y)dy + \sum_{i=1}^n\lambda_i p_i -\gamma \sum_{i=1}^n  \int_\mathcal{Y}\exp^\frac{g(y)+\lambda_i - c(y,z_i)}{\gamma} p_iq(y)dy,\right\}
%\end{align}
%where $c(y,z_i)$ is a cost function evaluated on the support $\mu$, $g$ is a dual potential and $\lambda$ is a dual variable. \\
%\dd{Here we duplicated the variable $p$ to execute the problem in a distributed manner.}\\
%We denote ${{{\mathtt{p}} = [p_1^T,\cdots,p_m^T]^T}}$, where $\forall ~i\in V$\dd{? (V wasn't defined so far)}, $p_i\in S_1(n)$\\

\subsection{Network constraints in the distributed barycenter problem}
%Consider continuous measures $\mu_1, .., \mu_n$ on $\mathcal{Y}$ with densities $q_i(y)$. Solving the Wasserstein barycenter problem we wish to approximate the true barycenter $\nu$ of these measures by its discrete approximation $\nu = \sum_{i=1}^np_i\delta(z_i)$.
%For $i=1,...,m$, consider measures $\mu_i$ supported on a metric spaces $\mathcal{Y}_i$ and having densities $q_i(y)$. 

We now describe the distributed optimization setting for solving problem~\eqref{w_barycenter}.
To do so, we rewrite the problem \eqref{w_barycenter} in an equivalent form
\begin{align}
\min_{\substack{p_1=\cdots=p_m \\ p_1,\dots,p_m \in S_1(n)}}  ~ %\W_{\gamma,\mu}(p) :=
\sum\limits_{i=1}^{m} \W_{\gamma,\mu_i}(p_i).
\label{eq:DecDistrPreProb}
\end{align}
We assume that each measure $\mu_i$ is held by an agent $i$ on a network and this agent can sample from this measure. We model such a network as a fixed \textit{connected undirected graph} \mbox{$\mathcal{G} = (V,E)$}, where 
$V$ is the set of $m$ nodes and $E$ is the set of edges. We assume that the graph $\mathcal{G}$ does not have self-loops. The network structure imposes information constraints, specifically, each node $i$ has access to $\mu_i$ only and can exchange information only with its immediate neighbors, i.e. nodes $j$ s.t. $(i,j)\in E$. 

We represent the communication constraints imposed by the network by introducing a single equality constraint instead of the constraints $p_1=\cdots=p_m$ in \eqref{eq:DecDistrPreProb}.
To do so, we define the Laplacian matrix 
$\bar W{\in \mathbb{R}^{m\times m}}$ of the graph $\mathcal{G}$ as
{\small
\begin{align*}
[\bar W]_{ij} = \begin{cases}
-1,  & \text{if } (i,j) \in E,\\
\text{deg}(i), &\text{if } i= j, \\
0,  & \text{otherwise,}
\end{cases}
\end{align*}}
\noindent\hspace{-2mm} where $\text{deg}(i)$ is the degree of the node $i$, i.e., the number of neighbors of the node. 
Finally,  we define the communication matrix (also referred to as an interaction matrix) 
by $W := \bar W \otimes I_n$, where $\otimes$ denotes the Kronecker product of matrices. 

%We define the degree condition number of $W$ as \mbox{$\chi(W) = {d_{\max}}/{d_{\min}} $}.

%Throughout the paper, {\it we assume that graph $\mathcal{G} = (V,E)$ is undirected and connected}.

Since the graph $\mathcal{G}$ is undirected and connected, the Laplacian matrix $\bar W$ is symmetric and positive semidefinite. Furthermore,
the vector $\boldsymbol{1}$ of all ones is the unique (up to a scaling factor) eigenvector associated with the zero eigenvalue. In respect that the matrix $W$ inherits the properties of $\bar W$, i.e., it is symmetric and positive, we conclude that
\[ W{\mathtt{p}} = 0  \;  \; \text{if and only if} \; \;  {p_1 = \cdots = p_m}, \]
where $\mathtt{p} = [p_1^T,\cdots,p_m^T]^T \in \mathbb{R}^{mn}$. 
Moreover, this identity holds for $\sqrt{W} := \sqrt{\bar W} \otimes I_n$, i.e. 
\[\sqrt{W}{\mathtt{p}} = 0  \; \;\text{if and only if}  \; \; {p_1 = \cdots = p_m}.\] 
%the following holds 
%\begin{itemize}
%    \setlength{\itemsep}{0pt}
%    \item $W{x} = 0$ if and only if {$x_1 = \cdots = x_m$}.
%    \item $\sqrt{W}{x} = 0$ if and only if {$x_1 = \cdots = x_m$}.
%\end{itemize}
% Note that the constraint set $ \{p_1,\dots, p_m \in S_1(n) \mid \sqrt{W} \mathtt{p}=0\}$
% is 
% the same as the set $\{p_1,\dots, p_m \in S_1(n) \mid p_1 = \cdots  = p_m\}$,
% since
% \mbox{$\ker (\sqrt{W}) = \spn(\boldsymbol{1})$} due to the connectivity of the graph $\mathcal{G}$. Here we defined stacked column vector $\mathtt{p} = [p_1^T,\cdots,p_m^T]^T \in \mathbb{R}^{mn}$. 
Using this fact, we equivalently rewrite problem~\eqref{w_barycenter} as the maximization problem with linear equality constraint
    \begin{align}\label{consensus_problem2}
    \max_{\substack{p_1,\dots, p_m \in S_1(n) \\ \sqrt{W} \mathtt{p}=0 }} ~ - \sum\limits_{i=1}^{m} \W_{\gamma, \mu_i}(p_i) .
    \end{align}
    
\subsection{Dual formulation of the barycenter problem}
    
%The Lagrangian dual problem for \eqref{consensus_problem2} is
%{\small
%\begin{align*}
%\min_{\lambda_1,\dots, \lambda_m \in \R^n} ~\max_{p_1,\dots, p_m \in S_1(n) } ~ \left\lbrace \sum\limits_{i=1}^{m} \langle \lambda_i, [\sqrt{W}\mathtt{p}]_i\rangle-\W_{\gamma, \mu_i}(p_i)\right\rbrace = \min_{\lambda_1,\dots, \lambda_m \in \R^n} \sum_{i=1}^{m} \W^*_{\gamma, \mu_i}([\sqrt{W}\vec{\lambda}]_i),
%\end{align*}}
%where $\vec{\lambda} = [\lambda_1^T,\cdots,\lambda_m^T]^T$ and $\W^*_{\gamma,\mu_i}(\lambda_i)$ is the Fenchel-Legendre transform of $\W_{\gamma,\mu_i}(p_i)$.
%Here we used that 
%\begin{align*}
%\mathcal{W}^*_{\gamma,\mu_i}([\sqrt{W}\mathtt{y}]_i) & = \max_{p_i \in S_1(n)}  \left\lbrace \left\langle [\sqrt{W}\mathtt{y}]_i,p_i \right\rangle -\mathcal{W}_{\gamma,q_i}(p_i) \right\rbrace,
%\end{align*}
%\noindent where $[\sqrt{W}y]_i$ is equivalent %form of $\sum_j^m \sqrt{W}_{ij}y_j$ and %$\sqrt{W}_{ij} = [\sqrt{\bar W}]_{ij} \otimes %I_n$. 
%By the chain rule, we can explicitly express 
%{\small\begin{align*}
%\nabla %\left(\W^*_{\gamma,\mu_i}\left([\sqrt{W}\vec{%\lambda}]_i\right) \right) & = %\sum_{j=1}^{m}\sqrt{W}_{ij} \nabla %\W_{\gamma,\mu_j}^*%([\sqrt{W}\vec{\lambda}]_i). 
%\end{align*}}

Given that problem \eqref{consensus_problem2} is an optimization problem with linear constraints, we introduce a vector of dual variables $\Blm = [\lambda_1^T,\cdots,\lambda_m^T]^T \in \R^{mn}$ for the constraints $\sqrt{W}\p=0$ in \eqref{consensus_problem2}.
Then, the Lagrangian dual problem for \eqref{consensus_problem2} is
\begin{align}\label{eq:DualPr}
\min_{\Blm \in \R^{mn}} ~\max_{p_1,\dots, p_m \in S_1(n) } ~ \left\lbrace \sum\limits_{i=1}^{m} \langle \lambda_i, [\sqrt{W}\mathtt{p}]_i\rangle-\W_{\gamma, \mu_i}(p_i)\right\rbrace = \min_{\Blm \in \R^{mn}} \sum_{i=1}^{m} \W^*_{\gamma, \mu_i}([\sqrt{W}\Blm]_i),
\end{align}
where $[\sqrt{W}\p]_i$ and $[\sqrt{W}\Blm]_i$ denote the $i$-th $n$-dimensional block of vectors $\sqrt{W}\p$ and $\sqrt{W}\Blm$ respectively, the equality $\sum\limits_{i=1}^{m} \langle \lambda_i, [\sqrt{W}\mathtt{p}]_i\rangle = \sum\limits_{i=1}^{m} \langle [\sqrt{W} \Blm]_i, p_i\rangle$ was used, and $\W^*_{\gamma,\mu_i}(\cdot)$ is the Fenchel-Legendre transform of $\W_{\gamma,\mu_i}(p_i)$. 
The following Lemma states that each $\W^*_{\gamma, \mu_i}(\cdot)$ is a smooth function with Lipschitz-continuous gradient and can be expressed as an expectation of a  function of additional random argument. 

\begin{lemma}\label{Lm:dual_obj_properties}
Given a positive Radon probability measure  $\mu \in \mathcal{M}_{+}^1(\mathcal{Y})$  with density $q(y)$  on a metric space $\mathcal{Y}$, the Fenchel-Legendre dual function for $\W_{\gamma,\mu}(p)$ has the following explicit form 
%{\small \begin{align*}
%\W_{\gamma,\mu}^*(\lambda) = \gamma \int_{\mathcal{Y}} q(y)\cu{\log}\left(\frac{1}{q(y)}\sum_{i=1}^n\exp\left(\frac{\lambda_i-c_i(y)}{\gamma} \right)\right) dy = \E_{Y\sim\mu}\gamma \cu{\log}\left(\frac{1}{q(y)}\sum_{i=1}^n\exp\left(\frac{\lambda_i-c_i(Y)}{\gamma} \right)\right).
%\end{align*}}
\begin{align*}
\W_{\gamma,\mu}^*(\blm) & = \E_{Y\sim\mu}\gamma \log\left(\frac{1}{q(Y)}\sum_{\ell=1}^n\exp\left(\frac{[\blm]_\ell-c_\ell(Y)}{\gamma} \right)\right),
\end{align*}
and its gradient is $1/\gamma$-Lipschitz continuous w.r.t. 2-norm with following $l$-th component 
 \begin{align*}
[\nabla \W_{\gamma,\mu}^*(\bar{\lambda})]_l = \E_{Y\sim\mu} \frac{\exp(([\bar{\lambda}]_l-c_l(Y))/\gamma)  }{\sum_{\ell=1}^n\exp(([\bar{\lambda}]_\ell-c_\ell(Y))/\gamma)}, ~l=1,\dots,n,
\end{align*}
where $Y\sim \mu$ means that random variable $Y$ is distributed according to measure $\mu$.
\end{lemma}

% \begin{align*}
% \mathcal{W}^*_{\gamma,\mu_i}(\blm_i) & = \max_{p_i \in S_1(n)}  \left\lbrace \left\langle \blm_i,p_i \right\rangle -\mathcal{W}_{\gamma,\cu{\mu_i}}(p_i) \right\rbrace,
% \end{align*}
% where $\blm_i := [\sqrt{W}\Blm]_i$, $i=1,...,m$.

Denote $\bar{\Blm} = \sqrt{W}\Blm = [[\sqrt{W}\Blm]^T_1, \dots, [\sqrt{W}\Blm]^T_m]^T=  [\blm^T_1, \dots, \blm^T_m]^T $ and  $\W^*_\gamma (\Blm)$ -- the dual objective in the r.h.s. of \eqref{eq:DualPr}. Then, by the chain rule, the $l$-th $n$-dimensional block of $\nabla \W^*_\gamma (\Blm)$ is
%by the Demyanov-Danskin's theorem, $ \nabla \W^*_{\gamma, \mu_i} (\blm_i) = p_i(\blm_i)$, where $p_i(\blm_i)$ is the unique solution to the inner maximization problem  
%\begin{align}\label{eq:primal-dual-var-connect}
%p_i(\blm_i) = \arg\max\limits_{p_i \in S_1(n)} \left\lbrace \left\langle \blm_i,p_i \right\rangle -\mathcal{W}_{\gamma,\cu{\mu_i}}(p_i) \right\rbrace.
%\end{align}
% By the chain rule,
% the gradient of the objective function $\W^*_\gamma (\sqrt{W}\Blm) = \sum\limits_{i=1}^{m} \W^*_{\gamma, \mu_i}([\sqrt{W}\Blm]_i)$ in the dual problem \eqref{eq:DualPr} is
% \begin{align}\label{eq:DualObjGrad-0}
% \nabla \W^*_\gamma (\sqrt{W}\Blm) =  \sum_{i=1}^{m} \nabla \W^*_{\gamma, \mu_i}([\sqrt{W}\Blm]_i) = \sqrt{W}\nabla \W^*_\gamma ({\BBlm}).
% \end{align}
% Moreover, its $i$-th $n$-dimensional block is
\begin{align}\label{eq:DualObjGrad}
\left[\nabla \W^*_\gamma (\Blm) \right]_l = 
\left[\nabla \sum_{i=1}^{m} \W^*_{\gamma, \mu_i}([\sqrt{W}\Blm]_i) \right]_l = \sum_{j=1}^{m}\sqrt{W}_{lj} \nabla \W_{\gamma,\mu_j}^*(\blm_j), \; l=1,...,m.
\end{align}

From Lemma~\ref{Lm:dual_obj_properties} and the expression  \eqref{eq:DualObjGrad} for the gradient of the dual objective, we can see that the \emph{dual} problem \eqref{eq:DualPr} is a smooth stochastic convex optimization problem. 
This is in contrast to \cite{lan2017communication}, where the primal problem is a stochastic optimization problem. Moreover, as opposed to the existing literature on stochastic convex optimization, we not only need to solve the dual problem but also need to reconstruct an approximate solution for the primal problem \eqref{consensus_problem2}, which is the barycenter.
%The difference with the existing literature on stochastic convex optimization is that we not only need to solve the dual problem but also need to reconstruct an approximate solution for the primal problem \eqref{consensus_problem2}, which is the barycenter. 
In order to do this, in the next section, we develop a novel accelerated primal-dual stochastic gradient method for a general smooth stochastic optimization problem, which is dual to some optimization problem with linear equality constraints. Further, in Section \ref{sec:GenAlgApplication}, we apply our general algorithm to the particular case of primal-dual pair of problems \eqref{consensus_problem2} and \eqref{eq:DualPr}.

%In the next problem we will consider a general problem of this type with an accelerated algorithm for its solution and convergence rate analysis. We will also show, how the primal variable can be reconstructed.

\section{General Primal-Dual Framework for Stochastic Optimization}\label{sec:stoch_prob}
In this section, we consider a general smooth stochastic convex optimization problem which is dual to some optimization problem with linear equality constraints. Extending our works \cite{dvurechensky2016stochastic,gasnikov2016efficient,chernov2016fast,dvurechensky2016primal-dual,dvurechensky2017adaptive,anikin2017dual,dvurechensky2018computational}, we develop a novel algorithm for its solution and reconstruction of the primal variable together with convergence rate analysis. 
%This algorithm and its analysis can be of independent interest. 
We underline that the material of this section is not standard. Unlike prior works, we consider the stochastic primal-dual pair of problems and one of our contributions consists in providing a primal-dual extension of the accelerated stochastic gradient method. We believe that our algorithm can be used for problems other than regularized Wasserstein barycenter problem and, thus, we, first, provide a general algorithm and, then, apply it to the barycenter problem. 
%The main difference of our setting from \cite{lan2017communication} is that in that paper the authors solve primal stochastic optimization problem, whereas in or case the \emph{dual} problem is a stochastic optimization problem.
We introduce new notation since this section is independent of the others and is focused on a general optimization problem.

\subsection{General setup and assumptions}
\label{S:stoch_prob_setup}
For any finite-dimensional real vector space $E$, we denote by $E^*$ its dual, by $\la \lambda, x \ra$ the value of a linear function $\lambda \in E^*$ at $x\in E$. Let $\|\cdot\|$ denote some norm on $E$ and $\|\cdot\|_{*}$ denote the norm on $E^*$ which is dual to $\|\cdot\|$, i.e.
$\|\lambda\|_{*} = \max \{ \la \lambda, x \ra :  \|x\| \leq 1 \}$.
%By $\partial f(x)$ we denote the subdifferential of a function $f(x)$ at a point $x$. 
% Let $E_1, E_2$ be two finite-dimensional real vector spaces. 
% For a linear operator $A:E_1 \to E_2$, we define its norm as follows
% $$
% \|A\|_{E_1 \to E_2} = \max_{x \in E_1,u \in E_2^*} \{\la u, A x \ra : \|x\|_{E_1} = 1, \|u\|_{E_2,*} = 1 \}.
% $$
For a linear operator $A:E_1 \to E_2$, we define the adjoint operator $A^T: E_2^* \to E_1^*$ in the following way $\la u, A x \ra = \la A^T u, x \ra, \quad \forall ~ u \in E_2^*, x \in E_1$.
We say that a function $f: E \to \R$ has a $L$-Lipschitz continuous gradient w.r.t. norm $\|\cdot\|_{*}$ if it is continuously differentiable and its gradient satisfies Lipschitz condition $
\|\nabla f(x) - \nabla f(y) \|_{*} \leq L \|x-y\|, \quad \forall ~x,y \in E$.
Note that, from this inequality, it follows that
\begin{equation}
f(y) \leq f(x) + \la \nabla f(x) , y-x \ra + \frac{L}{2} \|x-y\|^2, \quad \forall ~x,y \in E.
\label{eq:nfLipDef}
\end{equation}

The main problem, we consider in this section, is a  
where $Q$ is a simple closed convex set, $A: E \to H$ is given linear operator, $b \in H$ is given, $\Lambda = H^*$.
% The Lagrange dual problem to Problem $(P)$, written as a minimization problem is
% \begin{align}
% & (D) \quad \min_{\lambda \in \Lambda} \left\{   \la \lambda, b \ra +  + \max_{x\in Q} \left( -f(x) - \la A^T \lambda ,x \ra \right) \right\}. \notag
% %\label{eq:P_2}
% \end{align}
We define
\begin{equation}
\vp(\lambda) :=\la \lambda, b \ra +  \max_{x\in Q} \left( -f(x) - \la A^T \lambda  ,x \ra \right) = \la \lambda, b \ra + f^*(-A^T\lambda)
\label{eq:vp_def}
\end{equation}
and assume it to be smooth with $L$-Lipschitz continuous gradient. Here $f^*$ is the Fenchel-Legendre conjugate function for $f$. %Hence, $\nabla \vp(\lambda) = b - A x (\lambda)$, where $x (\lambda)$ is a solution of the problem
% Its gradient can be express as 
% \begin{equation}\label{eq:vp_grad}
% \nabla \vp(\lambda) = b - A x (\lambda),
% \end{equation}
% where $x (\lambda)$ is a solution of the problem
% \begin{equation}
% \max_{x\in Q} \left( -f(x) - \la A^T \lambda ,x \ra \right).
% \label{eq:inner}
% \end{equation}
%\pd{We need the following assumptions to prove Theorems 1 and 2. But we don't check them in the proof of Theorem 3. This should be done carefully. For example, here we assume that $\vp$ is a Fenchel conjugate of some function. But the dual objective $\W^*_{\gamma}$ in the barycenter problem is not dual to something, but is a  sum of duals.} 
We also assume that $f^*(-A^T\lambda)= \E_\xi F^*(-A^T\lambda,\xi)$, where $\xi$ is random vector and $F^*$ is the Fenchel-Legendre conjugate function to some function $F(x,\xi)$, i.e. it satisfies $F^*(-A^T\lambda,\xi) = \max\limits_{x\in Q}\{\la -A^T \lambda,x\ra - F(x,\xi) \}$. $F^*(\blm,\xi)$ is assumed to be smooth and, hence $\nabla_{\blm} F^*(\blm,\xi) = x(\blm,\xi)$, where $x(\blm,\xi)$ is the solution of the maximization problem
\begin{align*}
x(\blm,\xi) = \arg\max\limits_{x\in Q}\{\la \blm,x\ra - F(x,\xi)\}.
\end{align*}
%We also assume that $\nabla \vp(\lambda)$ is $L$-Lipschitz-continuous and 
Further, we assume that the dual problem $(D)$ can be accessed by a stochastic oracle $(\Phi(\lambda,\xi),~ \nabla \Phi(\lambda,\xi))$ with $\Phi(\lambda,\xi) = \la \lm, b\ra +F^*(-A^T\lambda,\xi)$ and $\nabla \Phi(\lambda,\xi) = b - A \nabla F^*(-A^T\lambda,\xi)$ satisfying
\begin{equation}
\E_\xi \Phi(\lambda,\xi) = \vp(\lambda), \quad \E_\xi \nabla \Phi(\lambda,\xi) = \nabla \vp(\lambda), \quad \E_\xi \|\nabla \Phi(\lambda,\xi) - \nabla \vp(\lambda) \|_2^2 \leq \sigma^2, \lambda \in H^*.
\end{equation}
Finally, we assume that dual problem $(D)$ has a solution $\lambda^*$ and there exists some $R >0$ such that $\|\lambda^*\|_{2} \leq R < +\infty$.

\subsection{An accelerated stochastic gradient method}
To solve the primal-dual pair of problems $(P)-(D)$, our first step, which we do in this subsection, is to introduce and analyse an accelerated stochastic gradient method (see Algorithm~\ref{Alg:ASGD}) for a general stochastic optimization problem and obtain some basic properties of the generated sequences, see Theorem~\ref{Th:ASGDConv}. In the next subsection, we apply it to the dual problem $(D)$. Algorithm~\ref{Alg:ASGD} is close in its form to the one in \cite{lan2012optimal}, but we use a different analysis extending \cite{dvurechensky2017adaptive} for the stochastic case.

To describe our algorithm, we introduce {\it proximal setup}, which is usually used in proximal gradient methods, see e.g. \cite{ben-tal2015lectures}. We choose some norm $\|\cdot\|$ on the space of vectors $\lambda$ and a {\it prox-function} $d(\lambda):\Lambda \rightarrow \R$ which is convex, continuous on $\Lambda$, continuously differentiable and $1$-strongly convex on $\Lambda_0 = \{\lm \in \Lambda: \partial d(\lm) \neq \emptyset\}$ with respect to $\|\cdot\|$, i.e., $\forall~ \lambda \in \Lambda, \zeta \in \Lambda^0$~ $d(\lm)-d(\zeta) -\la \nabla d(\zeta) ,\lm-\zeta \ra \geq \frac12\|\lm-\zeta\|^2$. Here $\partial d(\lm)$ is the subdifferential of $d$ and $\nabla d(x)$ is its subgradient.
%\begin{equation}
%d(y)-d(x) -\la \nabla d(x) ,y-x \ra \geq \frac12\|y-x\|^2.
%\label{eq:sc_def}
%\end{equation}
We define also the corresponding {\it Bregman divergence} $V[\zeta] (\lambda) := d(\lambda) - d(\zeta) - \la \nabla d(\zeta), \lambda - \zeta \ra$, $\lambda \in \Lambda, \zeta \in \Lambda^0$. 
It is easy to see that
\begin{equation}
V[\zeta] (\lambda) \geq \frac12\|\lambda - \zeta\|^2, \quad \forall~ \lambda \in \Lambda, \zeta \in \Lambda^0.
\label{eq:BFLowBound}
\end{equation}

\begin{algorithm}[h!]
\caption{Accelerated Stochastic Gradient Method (ASGD)}
\label{Alg:ASGD}
{\small
\begin{algorithmic}[1]
        \REQUIRE Starting point $\lambda_0 \in \Lambda$, prox-setup: $d(\lambda)$ -- $1$-strongly convex w.r.t. $\|\cdot\|$, the number of iterations $N$, Bregman divergence $V[\zeta] (\lambda) := d(\lambda) - d(\zeta) - \la \nabla d(\zeta), \lambda - \zeta \ra$, $\lambda \in \Lambda, \zeta \in \Lambda^0$.
        %\STATE Set $k=0$, 
        \STATE $C_0=\alpha_0=0$, $\eta_0=\zeta_0=\lambda_0$.
        %\REPEAT
        \FOR{$k = 0,\dots,N-1 $}
            %\STATE Set $M_k=L_k/2$.
            %\REPEAT
                %\STATE Set $M_k=2M_k$, 
                \STATE Find $\alpha_{k+1}$ as the largest root of the equation
                \begin{equation}
                C_{k+1}:=C_k+\alpha_{k+1} = 2L\alpha_{k+1}^2.\label{eq:alpQuadEq}
                \end{equation}
 %               \STATE 
%                \[M_{k+1} = \max \left\{1, ~\ceil[\Bigg]{\frac{\sigma^2C_{k+1}}{L\alpha_{k+1}\e}}\right\}\]
                %and 
                %\begin{equation}
                %C_{k+1}=C_k+\alpha_{k+1}
                %\label{eq:CkDef}
                %\end{equation}
                \STATE %Calculate 
                \begin{equation}
                \lambda_{k+1} = \frac{\alpha_{k+1}\zeta_k + C_k \eta_k}{C_{k+1}} .
                \label{eq:lambdakp1Def}
                \end{equation}
                \STATE %Calculate 
                \begin{equation}
                \zeta_{k+1}=\arg \min_{\lambda \in \Lambda} \{V[\zeta_{k}](\lambda) + \alpha_{k+1}(\Phi(\lambda_{k+1}, \xi_{k+1}) + \langle \nabla \Phi(\lambda_{k+1},\xi_{k+1}), \lambda - \lambda_{k+1} \rangle) \}.\label{eq:zetakp1Def}
                \end{equation}
                \STATE %Calculate
                \begin{equation}
                \eta_{k+1} = \frac{\alpha_{k+1}\zeta_{k+1} + C_k \eta_k}{C_{k+1}}.\label{eq:etakp1Def}
                \end{equation}
                \ENDFOR
            %\UNTIL{
            %\begin{equation}
            %\vp(\eta_{k+1}) \leq \vp(\lambda_{k+1}) + \la \nabla \vp(\lambda_{k+1}) ,\eta_{k+1} - \lambda_{k+1} \ra  +\frac{M_k}{2}\|\eta_{k+1} - \lambda_{k+1}\|^2.
            %\label{eq:lipConstCheck}
            %\end{equation}}
            %\STATE Set $L_{k+1}=M_k/2$, $k=k+1$.
        %\UNTIL{Option 1: $k = k_{\max}$. \\ 
                    %Option 2: $R^2/C_k \leq \e$. \\
                    %Option 3:  
                    %$$\vp(\eta_k) - \min_{\lambda \in \Lambda: V[\zeta_0](\lambda) \leq R^2} \left\{ \sum_{i=0}^k \frac{\alpha_{i}}{C_k} \left( \Phi(\lambda_{i}) + \la \nabla \Phi (\lambda_{i}), \lambda - \lambda_{i}\ra \right) \right\} \leq \e.$$
                    %Here $R$ is such that $ V[\zeta_0](\lambda_*) \leq R^2$ and $\e$ is the desired accuracy.}
        \ENSURE The point $\eta_{N}$.    
\end{algorithmic}
}
\end{algorithm}

\begin{theorem}
\label{Th:ASGDConv}
Let the sequences $\{\lambda_N, \eta_N, \zeta_N, \alpha_N, C_N \}$, $N>0 $ be generated by Algorithm \ref{Alg:ASGD}. Then,  for all $N > 0$, it holds that
{\small \begin{align} 
    C_N\vp(\eta_N) 
    &\leq \min_{\lambda \in \Lambda} \left\{ \sum_{k=0}^N \alpha_{k} \left( \vp(\lambda_{k}) + \la \nabla  \Phi(\lambda_{k},\xi_k), \lambda - \lambda_{k}\ra \right) + V[\zeta_0](\lambda) \right\} \notag \\ 
 &\hspace{-3em}+ \sum_{k=0}^{N-1}C_{k}\la \nabla \Phi(\lambda_{k+1},\xi_{k+1}) - \nabla \vp(\lambda_{k+1}), \eta_{k} - \lambda_{k+1}\ra + \sum_{k=0}^{N}\frac{C_{k}}{2L}\| \nabla  \Phi(\lambda_{k},\xi_k)- \nabla \vp(\lambda_{k})\|_{*}^2.\label{eq:ASGDConv}
    \end{align}}
\end{theorem}
%\pd{TODO:Add proofs.}

%\begin{lemma}
%\label{Lm:CkGrowth}
%Let the sequence $\{C_k \}$, $k\geq 0$ be generated by Algorithm \ref{Alg:ASGD}. Then,  for all $k \geq 1$ it holds that
%\begin{equation}
%C_k \geq \frac{(k+1)^2}{8L},
%\label{eq:CkGrowth}
%\end{equation}	
%where $L$ is the Lipschitz constant for the gradient of $\vp$.
%\end{lemma}
%\begin{proof}For $k = 1$ \eqref{eq:CkGrowth} holds 
%	\begin{equation*}
%		C_1 = \alpha_1 \geq \frac{1}{2L}.
%	\end{equation*}
 %   Further, based on induction one can show that it holds for all $k$.
%\end{proof}

\subsection{Accelerated primal-dual stochastic gradient method}
In this subsection, we develop an accelerated algorithm for the primal-dual pair of problems $(P)-(D)$. The idea is to apply the algorithm of the previous subsection to the dual problem $(D)$, endow it with a step in the primal space and, using the result of Theorem \ref{Th:ASGDConv}, show that the new algorithm allows to approximate also the solution to the primal problem. Since the feasible set in the problem $(D)$ is unbounded, we choose the Euclidean proximal setup in $H^*$ and denote the standard Euclidean norm by $\|\cdot\|_2$. We use Euclidean proximal setup with the prox-function  $d(\lambda) = \frac{1}{2}\|\lambda\|^2_2$ and the Bregman divergence  $V[\zeta](\lambda) = \frac{1}{2}\|\lambda - \zeta\|^2_2$.

Note that, in this case, the dual norm is also Euclidean and the step 5 of the algorithm simplifies. We additionally assume that the variance of the stochastic approximation $\nabla \Phi(\lambda,\xi)$ for the gradient of $\vp$ can be controlled and made as small as we desire. This can be done, for example by mini-batching the stochastic approximation. Finally, since $\nabla  \Phi(\lambda,\xi) = b - A \nabla F^*(-A^T\lambda,\xi) = b - A x(-A^T\lambda,\xi)$, on each iteration, to find $\nabla  \Phi(\lambda_,\xi)$ we find the vector $x(-A^T\lambda,\xi)$ and use it for the primal iterates.

\begin{algorithm}[h!]
\caption{Accelerated Primal-Dual Stochastic Gradient Method (APDSGD)}
\label{Alg:APDSGD}
{\small
\begin{algorithmic}[1]
   \REQUIRE starting point $\lambda_0 = 0$,  the number of iterations $N$.
   %, accuracy $\tilde{\e}_f,\tilde{\e}_{eq},\tilde{\e}_{in} > 0$.
   %\STATE Set $k=0$, 
   \STATE $C_0=\alpha_0=0$, $\eta_0=\zeta_0=\lambda_0= \hat{x}_0= 0$.
   %\REPEAT
            %\STATE Set $M_k=L_k/2$.
            %\REPEAT
                %\STATE Set $M_k=2M_k$, 
                \FOR{$k=0,\dots, N-1$}
                \STATE Find $\alpha_{k+1}$ as the largest root of the equation \begin{equation}
                C_{k+1}:=C_k+\alpha_{k+1} = 2L\alpha_{k+1}^2.
                \label{eq:PDalpQuadEq}
                \end{equation}
               % \STATE 
                %\begin{equation}\label{M}
                %M_{k+1} = \max \left\{1, ~\ceil[\Bigg]{\frac{\sigma^2C_{k+1}}{L\alpha_{k+1}\e}}\right\}
                %\end{equation}
                %and 
                %\begin{equation}
                %C_{k+1}=C_k+\alpha_{k+1}
                %\label{eq:CkDef}
                %\end{equation}
                \STATE
                \begin{equation}
                \lambda_{k+1} = \frac{\alpha_{k+1}\zeta_k + C_k \eta_k}{C_{k+1}} .
                \label{eq:PDlambdakp1Def}
                \end{equation}
                \STATE 
                \begin{equation}\zeta_{k+1}= \zeta_{k} - \alpha_{k+1} \nabla  \Phi(\lambda_{k+1},\xi_{k+1}).
                \label{eq:PDzetakp1Def}
                \end{equation}
                \STATE \begin{equation}
                \eta_{k+1} =\frac{\alpha_{k+1}\zeta_{k+1} + C_k \eta_k}{C_{k+1}}.
                \label{eq:PDetakp1Def}
                \end{equation}
            %\UNTIL{
            %    \begin{equation}
            %        \vp(\eta_{k+1}) \leq \vp(\lambda_{k+1}) + \la \nabla \vp(\lambda_{k+1}) ,\eta_{k+1} - \lambda_{k+1} \ra  +\frac{M_k}{2}\|\eta_{k+1} - \lambda_{k+1}\|_2^2.
            %        \label{eq:PDlipConstCheck}
            %    \end{equation}}
            \STATE Set
                    \begin{equation}
                    \hat{x}_{k+1} = \frac{1}{C_{k+1}}\sum_{i=0}^{k+1} \alpha_i x(-A^T\lambda_i,\xi_i) = \frac{\alpha_{k+1}x(-A^T\lambda_{k+1},\xi_{k+1})+C_k\hat{x}_{k}}{C_{k+1}}.\notag
                \label{eq:xhat_def}
                \end{equation}
                \ENDFOR
            %\STATE Set $L_{k+1}=M_k/2$, $k=k+1$.
  %\UNTIL{$|f(\hat{x}_{k+1})+\vp(\eta_{k+1})| \leq \tilde{\e}_f$, $\|A_1\hat{x}_{k+1}-b_1\|_{2} \leq \tilde{\e}_{eq}$, $\rho(A_2\hat{x}_{k+1}-b_2,-K) \leq \tilde{\e}_{in}$.}
    \ENSURE The points $\hat{x}_{N}$, $\eta_{N}$.    
\end{algorithmic}}
\end{algorithm}

\begin{theorem}
\label{Th:stoch_err}
Let $\vp$ have $L$-Lipschitz continuous gradient w.r.t. 2-norm and $\|\lambda^*\|_2 \leq R$, where $\lambda^*$ is a solution of dual problem $(D)$. Assume that at each iteration of Algorithm \ref{Alg:APDSGD}, the stochastic approximation $\nabla \Phi(\lambda_k,\xi_k)$ of the gradient is chosen in such a way that $\E_\xi \|\nabla \Phi(\lambda_k,\xi_k) - \nabla \vp(\lambda_k) \|_2^2 \leq \frac{\e L\alpha_k}{C_k}$.
Then, for any $\varepsilon > 0$ and $N \geq 0$, the output $\hat{x}_N$ generated by the  Algorithm \ref{Alg:APDSGD} satisfies  
\begin{align}\label{eq:APDSGDRate}
%\mathbb{E}_{\xi_1,\dots, \xi_N}[\vp(\eta_N)+f(\hat{x}_N)] \leq \frac{16LR^2}{N^2} + \frac{\e}{2}~ ~ ~ ~ \text{and} ~ ~ ~ ~ \mathbb{E}_{\xi_1,\dots, \xi_N}\|A\hat{x}_N-b\|_2 \leq \frac{16LR}{N^2} + \frac{\e}{2R}.
f(\mathbb{E}\hat{x}_N)-f^* \leq \frac{16LR^2}{N^2} + \frac{\e}{2}~ ~ ~ ~ \text{and} ~ ~ ~ ~ \|A\mathbb{E}\hat{x}_N-b\|_2 \leq \frac{16LR}{N^2} + \frac{\e}{2R},
\end{align}
where the expectation is taken w.r.t. all the randomness $\xi_1,\dots, \xi_N$.
% Let $f$ is $\gamma$-strongly convex and $\|\lambda^*\|_2 \leq R$, where $\lambda^*$ is a solution of dual problem (D), then after $N\geq \sqrt{\frac{16LR^2}{\varepsilon}}$ iterations of Algorithm \ref{Alg:APDSGD} the following 
% \begin{align*}
% \mathbb{E}_{\xi_1,\dots, \xi_N}[\vp(\eta_N)+f(\hat{x}_N)] \leq 2\varepsilon ~ ~ ~ ~ \text{and} ~ ~ ~ ~ \mathbb{E}_{\xi_1,\dots, \xi_N}\|A\hat{x}_N-b\|_2 \leq 2\varepsilon/R
% \end{align*}
% holds with stochastic error $\sigma_k^2 = \frac{\e L\alpha_k}{C_k}$, where  ~$
% \mathbb{E}_{\xi}\left[\left\| \nabla \Phi(\lambda_{k},\xi_k)- \nabla \vp(\lambda_{k})\right \|_2^2 \right] \leq \sigma^2_k.$
\end{theorem}
%\dd{\subsection{Mini-batch stochastic gradient}}

In step 7 of Algorithm \ref{Alg:APDSGD} we can use a batch of size $M$ and $\frac{1}{M}\sum_{r=1}^M x(\lambda_{k+1},\xi_{k+1}^r)$ to update $\hat{x}_{k+1}$. Then, under reasonable assumptions, $\hat{x}_N$ concentrates around $\mathbb{E}\hat{x}_N$ \cite{guigues2017non-asymptotic}  and, if $f$ is Lipschitz-continuous, we obtain that \eqref{eq:APDSGDRate} holds with large probability with $\hat{x}_N$ instead of $\mathbb{E}\hat{x}_N$.

\section{Solving the Barycenter Problem}
\label{sec:GenAlgApplication}
In this section, we apply the general algorithm APDSGD from the previous section to solve the primal-dual pair of problems \eqref{consensus_problem2}-\eqref{eq:DualPr} and approximate the regularized Wasserstein barycenter which is a solution to \eqref{consensus_problem2}. 
First, in Lemma \ref{Lm:dual_obj_properties2}, we make a number of technical steps to take care of the assumptions of Theorem \ref{Th:stoch_err}. We estimate the Lipschitz constant of the dual objective's gradient in \eqref{eq:DualPr}, introduce mini-batch stochastic approximation for the gradient of the dual objective and estimate its variance. 
Then, we introduce a change of dual variable so that a gradient-type step for the dual objective, e.g., the step 5 of Algorithm~\ref{Alg:APDSGD}, becomes feasible for the decentralized distributed setting. 
Then, for simplicity, we consider a non-accelerated algorithm for regularized Wasserstein barycenter problem to illustrate the combination of gradient methods, a stochastic approximation of the gradient and decentralized distributed computations.
Finally, we present our accelerated algorithm for regularized Wasserstein barycenter problem with its complexity analysis.

%Next lemma estimates the Lipschitz constant of the dual objective in \eqref{eq:DualPr}, estimates mini-batch stochastic approximation for the gradient of the dual objective and estimates its variance. 

\begin{lemma}\label{Lm:dual_obj_properties2}
%\footnote{The scaling factor $m$ in the variance is due to initial assumption on all weights corresponding to Wasserstein distances in the sum are equal to 1, that increases the variance in $m$ times.}
The gradient of the dual objective function $\W_{\gamma}^*(\Blm)$ in the dual problem \eqref{eq:DualPr} is $\lambda_{\max}(W) /\gamma$-Lipschitz continuous w.r.t. 2-norm. If its stochastic approximation is defined as
\begin{align}
[\widetilde{\nabla} \W_{\gamma}^*(\Blm)]_i            &=  \sum_{j=1}^{m}\sqrt{W}_{ij} \widetilde{\nabla} \W_{\gamma,\mu_j}^*(\blm_j), \; i=1,...,m, \quad \text{with}\label{eq:tnW_def} \\
\widetilde{\nabla} \W_{\gamma,\mu_j}^*(\blm_j) &= \frac{1}{M}\sum_{r=1}^{M} p_j(\blm_j, Y^j_r), \; j=1,...,m, \quad \text{and}\label{eq:StochGradDef0} \\
[p_j(\blm_j, Y^j_r)]_l &=
\frac{\exp(([\bar{\lambda}_j]_l-c_l(Y_r^j))/\gamma)  }{\sum_{\ell=1}^n\exp(([\bar{\lambda}_j]_{\ell}-c_\ell(Y_r^j))/\gamma)}, \; j=1,...,m, \; l=1,...,n, \; r=1,...,M\label{eq:StochGradDef}
\end{align}    
% \begin{align}
% [\widetilde{\nabla} \W_{\gamma}^*(\Blm)]_i &=  \sum_{j=1}^{m}\sqrt{W}_{ij} \widetilde{\nabla} \W_{\gamma,\mu_j}^*(\blm_j), \; i=1,...,m, \; \text{with}\notag \\
% [\widetilde{\nabla} \W_{\gamma,\mu_j}^*(\blm_j)]_l &= \frac{1}{M}\sum_{r=1}^{M} [p_j(\blm_j)]_l:=
% \frac{\exp(([\bar{\lambda}_j]_l-c_l(Y_r^j))/\gamma)  }{\sum_{\ell=1}^n\exp(([\bar{\lambda}_j]_{\ell}-c_\ell(Y_r^j))/\gamma)}, \; j=1,...,m, \; l=1,...,n, \label{eq:StochGradDef}
% \end{align}
where $M$ is the batch size, %$\blm_j := [\sqrt{W}\Blm]_j$, $j=1,...,m$, 
$Y_1^j,...,Y_r^j$ is a sample from the measures $\mu_j$, $j=1,...,m$. Then  
\begin{align*}\label{eq:variance}
&\E \widetilde{\nabla} \W_{\gamma}^*(\Blm) = \nabla \W_{\gamma}^*(\Blm) \quad \mbox{and} \\
&\E \|\widetilde{\nabla} \W_{\gamma}^*(\Blm) - \nabla \W_{\gamma}^*(\Blm)\|_2^2 \leq \frac{\lambda_{\max}(W)m}{M}, \; \Blm \in \R^{mn}, 
\end{align*}
where the expectation is taken w.r.t. all samples $(Y_1^j,\dots Y_M^j )$ from measure $\mu_j$, $j=1, \dots, m.$
\end{lemma}

Let us consider a simple stochastic gradient step for the particular dual problem \eqref{eq:DualPr}. Note that the step 5 of Algorithm \ref{Alg:APDSGD} has the same form. Using \eqref{eq:DualObjGrad}, the stochastic gradient step 
$\Blm_{k+1}= \Blm_{k}- \frac{1}{L} \widetilde{\nabla} \W_{\gamma}^*(\Blm_{k}) $
can be written block-wise as 
\begin{align*}
[\Blm_{k+1}]_i= [\Blm_{k}]_i- \frac{1}{L} \sum_{j=1}^{m}\sqrt{W}_{ij} \widetilde{\nabla} \W_{\gamma,\mu_j}^*([\sqrt{W}\Blm_{k}]_j), \quad \mbox{for each agent $i=1,...,m$},
\end{align*}
where $\widetilde{\nabla} \W_{\gamma,\mu_j}^*(\cdot)$ is defined in \eqref{eq:StochGradDef0} with the batch size $M=1$, and $L=\lambda_{\max}(W) /\gamma$.
%where $\tilde{\nabla} \W_{\gamma,\mu_j}^*(\cdot)$ means a stochastic approximation of the gradient of $\W_{\gamma,\mu_j}^*(\cdot)$.
Unfortunately, this update can not be made in the decentralized setting since the sparsity pattern of $\sqrt{W}_{ij}$ can be different from $W_{ij}$ and this will require some agents to get information not only from their neighbors. To overcome this obstacle, we change the variable and denote $\bar{\Blm} = \sqrt{W}\Blm$. 
%, $\bar{\Beta} = \sqrt{W}\Beta$, $\bar{\Bzeta} = \sqrt{W}\Bzeta$.
Then the gradient step becomes
\begin{align*}
[\BBlm_{k+1}]_i= [\BBlm_{k}]_i - \frac{1}{L}\sum_{j=1}^{m}W_{ij} \widetilde{\nabla} \W_{\gamma,\mu_j}^*([\BBlm_{k}]_j), \quad \mbox{for each agent $i=1,...,m.$}
\end{align*}
Algorithm \ref{alg:non-accler-main} presents a non-accelerated primal-dual stochastic gradient method, combining distributed updates and stochastic gradient step described above. This algorithm solves the primal-dual pair of problems \eqref{consensus_problem2}-\eqref{eq:DualPr} and approximates the regularized Wasserstein barycenter which is a solution to \eqref{consensus_problem2}. 
The algorithm has a loop, indexed by iteration number $k$ and the index $i$ corresponds to the agent's number. At each iteration $k$ of the algorithm, each agent $i$ samples from the measure $\mu_i$ and forms a stochastic approximation of the gradient of $\W_{\gamma,\mu_i}(\cdot)$. Then each agent shares this vector with its neighbors. After that, each agent calculates a step direction based on its information and information gathered from the neighbors. Note that the matrix $W$ provides communications only between neighboring nodes and step 6 requires only local information.

\begin{algorithm}[ht]
    \caption{Non-accelerated Distributed Computation of Wasserstein barycenter}
    \label{alg:non-accler-main}
    {
    \begin{algorithmic}[1]
    \REQUIRE Each agent $i\in V$ is assigned its measure $\mu_i$.
        \STATE All agents set $ [\BBlm_{0}]_i = \boldsymbol{0} \in \mathbb{R}^n$, $[\hat{\p}_0]_i= \boldsymbol{0} \in \R^n$, and $N$. 
%        \STATE{Set $K = \exp(-M/\gamma)$}
        \STATE{For each agent $i \in V$:}
        \FOR{ $k=0,\dots,N-1$ }
        \STATE{Sample $Y^i$  from the measure $\mu_i$ and set $\widetilde{\nabla} \W_{\gamma,\mu_i}^*([\BBlm_{k}]_i)$} as defined in \eqref{eq:StochGradDef0} with $M=1$.
%         \STATE  for all $l=1,\dots, n$ $$%[s_i]_l
%        \dd{[\tilde{\nabla} \W_{\gamma,\mu_i}^*(\tilde \lambda)]_l} = \frac{1}{\cu{M_{k+1}}}\sum_{r=1}^\cu{M_{k+1}}\frac{\exp(([\tilde\lambda_{k+1}^i]_l-c(y^i_r,z_l))/\gamma)  }{\sum_{t=1}^n\exp(([\tilde\lambda_{k+1}^i]_t-c(y^i_r,z_t))/\gamma)}$$
        \STATE{Share %$s_i$ 
        $\widetilde{\nabla} \W_{\gamma,\mu_i}^*([\BBlm_{k}]_i)$ with $\{j \mid (i,j) \in E \}$}
        %\STATE Calculate$$\tilde \eta_{k+1}^i  = \tilde \lambda_{k+1}^i - \frac{1}{L } \sum_{j=1}^{m} W_{ij} s_j$$
        \STATE  $$ [\BBlm_{k+1}]_i = [\BBlm_{k}]_i-\frac{1}{L} \sum_{j=1}^{m} W_{ij} \widetilde{\nabla} \W_{\gamma,\mu_j}^*([\BBlm_{k}]_j).$$
   \STATE $[\hat{\p}_{k+1}]_i  = [\hat{\p}_{k}]_i + \frac{1}{N} p_i([\BBlm_{k+1}]_i, Y^i),$
   where $p_i(\cdot,\cdot)$ is defined in \eqref{eq:StochGradDef}. 
   \ENDFOR
            %\STATE Set $L_{k+1}=M_k/2$, $k=k+1$.
  %\UNTIL{$|f(\hat{x}_{k+1})+\vp(\eta_{k+1})| \leq \tilde{\e}_f$, $\|A_1\hat{x}_{k+1}-b_1\|_{2} \leq \tilde{\e}_{eq}$, $\rho(A_2\hat{x}_{k+1}-b_2,-K) \leq \tilde{\e}_{in}$.}
    \ENSURE $\hat{\p}_{N}$.%, ~$\hat{y}^i_N = \tilde\eta_N^i $.   
        %\ENSURE $( y_N^*)_i = \tilde w^i_N $, ~
        % $(p^*_N)_i = \sum_{k=0}^{N-1} \frac{(k+2)}{N(N+3)}p^*_i(\tilde y^i_{k+1})$, $\forall i\in V$       
        
        %\STATE{Set $(\tilde y_N^*)_i = \frac{1}{N (N+3)}\left(\sum_{k=1}^{N-1}\tilde y^i_{k} + (N+1)^2\tilde y^i_{N} \right) $.}
        %\STATE{Set $(p^*_N)_i = \sum_{k=0}^{N-1} \frac{2(k+2)}{N(N+3)}p^*_i(\tilde y^i_{k+1})$.}
    \end{algorithmic}}
\end{algorithm}

Finally, we apply accelerated primal-dual stochastic gradient method (APDSGD) from the previous section to solve the primal-dual pair of problems \eqref{consensus_problem2}-\eqref{eq:DualPr} and calculate the regularized Wasserstein barycenter. As above, we introduce the change of dual variables $\bar{\Blm} = \sqrt{W}\Blm$, $\bar{\Beta} = \sqrt{W}\Beta$, $\bar{\Bzeta} = \sqrt{W}\Bzeta$, which makes the step 5 of Algorithm \ref{Alg:APDSGD} feasible for the decentralized distributed setting. The result is Algorithm \ref{alg:main}.
At each iteration $k$ each agent $i$ generates a sample of size $M_{k}$ from measure $\mu_i$, forms a stochastic approximation of the gradient of $\W_{\gamma,\mu_i}(\cdot)$ according to \eqref{eq:StochGradDef0} and shares it with the neighbors. 
The mini-batch size $M_k$ is chosen such that $M_k \geq \frac{m\gamma C_{k}}{\alpha_k\e}$, which, by Lemma \ref{Lm:dual_obj_properties2}, means that $\E \|\widetilde{\nabla} \W^*_{\gamma}(\Blm) - \nabla \W^*_{\gamma}(\Blm) \|_2^2 \leq \frac{\e L\alpha_k}{C_k}$ and the assumptions of Theorem \ref{Th:stoch_err} hold. 
% Based on this lemma, we see that if, on each iteration of Algorithm \ref{Alg:APDSGD}, the mini-batch size $M_k$ satisfies $M_k \geq \frac{\lambda_{\max}(W)m\gamma C_{k}}{\alpha_k\e}$ with $L= \frac{\lm_{\max}(W)}{\gamma}$, the assumptions of Theorem \ref{Th:stoch_err} hold. 

\begin{algorithm}[ht]
    \caption{Accelerated Distributed Computation of Wasserstein barycenter}
    \label{alg:main}
    {
    \begin{algorithmic}[1]
    \REQUIRE Each agent $i\in V$ is assigned its measure $\mu_i$.
        \STATE All agents $i \in V$ set $[\BBeta_{0}]_i = [\BBzeta_{0}]_i = [\BBlm_{0}]_i = [\hat{\p}_0]_i = \boldsymbol{0} \in \mathbb{R}^n$, $C_0=\alpha_0=0$ and $N$. 
%        \STATE{Set $K = \exp(-M/\gamma)$}
        \STATE{For each agent $i \in V$:}
        \FOR{ $k=0,\dots,N-1$ }
        \STATE Find $\alpha_{k+1}$ as the largest root of the equation $C_{k+1}:=C_k+\alpha_{k+1} = \frac{2\lm_{\max}(W)\alpha_{k+1}^2}{\gamma}.$
        \STATE $$M_{k+1} = \max\left\{1,~\ceil[\Big]{\frac{m\gamma C_{k+1}}{\alpha_{k+1}\e}} \right\}.$$
       % \STATE Calculate $$\tau_k=\alpha_{k+1}/C_{k+1}.$$
        %\STATE Calculate $$\tilde \lambda_{k+1}^i = \tau_k\tilde \zeta_k^i + (1-\tau_k) \tilde \eta_k^i.$$
        \STATE $$[\BBlm_{k+1}]_i = \frac{\alpha_{k+1}[\BBzeta_k]_i + C_k [\BBeta_k]_i}{C_{k+1}}.$$
        
        \STATE{Generate $M_{k+1}$ samples $\{Y_r^i\}_{r=1}^{M_{k+1}}$  from the measure $\mu_i$ and set $\widetilde{\nabla} \W_{\gamma,\mu_i}^*([\BBlm_{k}]_i)$ as in \eqref{eq:StochGradDef0} with $M=M_k$.}
%         \STATE  for all $l=1,\dots, n$ $$%[s_i]_l
%        \dd{[\tilde{\nabla} \W_{\gamma,\mu_i}^*(\tilde \lambda)]_l} = \frac{1}{\cu{M_{k+1}}}\sum_{r=1}^\cu{M_{k+1}}\frac{\exp(([\tilde\lambda_{k+1}^i]_l-c(y^i_r,z_l))/\gamma)  }{\sum_{t=1}^n\exp(([\tilde\lambda_{k+1}^i]_t-c(y^i_r,z_t))/\gamma)}$$
        \STATE{Share %$s_i$ 
        $\widetilde{\nabla} \W_{\gamma,\mu_i}^*([\BBlm_{k+1}]_i)$ with $\{j \mid (i,j) \in E \}$.}
        %\STATE Calculate$$\tilde \eta_{k+1}^i  = \tilde \lambda_{k+1}^i - \frac{1}{L } \sum_{j=1}^{m} W_{ij} s_j$$
        \STATE  $$ [\BBzeta_{k+1}]_i = [\BBzeta_{k}]_i-\alpha_{k+1} \sum_{j=1}^{m} W_{ij} \widetilde{\nabla} \W_{\gamma,\mu_j}^*([\BBlm_{k+1}]_j).$$
        \STATE  $$[\BBeta_{k+1}]_i = \frac{\alpha_{k+1}[\BBzeta_{k+1}]_i + C_k [\BBeta_{k+1}]_i }{C_{k+1}}.$$
        %\STATE{Share $p^*_i(\tilde y^i_{k+1})$ with $\{j \mid (i,j) \in E \}$.}
   \STATE $$[\hat{\p}_{k+1}]_i  = \frac{1}{C_{k+1}}\sum_{i=0}^{k+1} \alpha_i p_i([\BBlm_{k+1}]_i, Y^i_1) = \frac{\alpha_{k+1}p_i([\BBlm_{k+1}]_i, Y^i_1) + C_k[\hat{\p}_{k}]_i}{C_{k+1}},$$ 
   where $p_i(\cdot,\cdot)$ is defined in \eqref{eq:StochGradDef}. \footnotemark 
   %\pd{It is not clear, which $r$ to take here. Any? Or maybe weshould batch also $p_i(\cdot,\cdot)$? I think we should average them.}
   \ENDFOR
            %\STATE Set $L_{k+1}=M_k/2$, $k=k+1$.
  %\UNTIL{$|f(\hat{x}_{k+1})+\vp(\eta_{k+1})| \leq \tilde{\e}_f$, $\|A_1\hat{x}_{k+1}-b_1\|_{2} \leq \tilde{\e}_{eq}$, $\rho(A_2\hat{x}_{k+1}-b_2,-K) \leq \tilde{\e}_{in}$.}
    \ENSURE $\hat{\p}_{N}$.%, ~$\hat{y}^i_N = \tilde\eta_N^i $.   
        %\ENSURE $( y_N^*)_i = \tilde w^i_N $, ~
        % $(p^*_N)_i = \sum_{k=0}^{N-1} \frac{(k+2)}{N(N+3)}p^*_i(\tilde y^i_{k+1})$, $\forall i\in V$       
        
        %\STATE{Set $(\tilde y_N^*)_i = \frac{1}{N (N+3)}\left(\sum_{k=1}^{N-1}\tilde y^i_{k} + (N+1)^2\tilde y^i_{N} \right) $.}
        %\STATE{Set $(p^*_N)_i = \sum_{k=0}^{N-1} \frac{2(k+2)}{N(N+3)}p^*_i(\tilde y^i_{k+1})$.}
    \end{algorithmic}}
\end{algorithm}
%\pd{TODO: define $p_i([\BBlm_{k+1}]_i)$ when the theorem 2 will be proved.}
\footnotetext{Note that we can use also $\frac{1}{M_{k+1}}\sum_{r=1}^{M_{k+1}} p_i([\BBlm_{k+1}]_i, Y^i_r)$ instead of $p_i([\BBlm_{k+1}]_i, Y^i_1)$. This does not change the statement of Theorem \ref{Th:WBCompl}, but reduces the variance of $\hat{\p}_{N}$ in practice. Thus, in the experiments, we use this estimator for the primal variable. Moreover, under mild assumptions, we can obtain high-probability analogue to inequalities \eqref{eq:WBerr}.}

\begin{theorem}\label{Th:WBCompl}
Let the assumptions of Section \ref{sec:problem} hold and $R$ be such that $\|\Blm^*\|_2 \leq R$. Then Algorithm \ref{alg:main} after $N=\sqrt{32\lambda_{\max}(W)R^2/(\e\gamma)}$ iterations returns an approximation $\hat{\p}_{N}$ for the barycenter, which satisfies
\begin{equation}\label{eq:WBerr}
\sum\limits_{i=1}^{m} \W_{\gamma,\mu_i}(\E[\hat{\p}_{N}]_i)-\sum\limits_{i=1}^{m} \W_{\gamma,\mu_i}([\p^*]_i) \leq \e, \quad \|\sqrt{W}\E\hat{\p}_{N}\|_2 \leq \e/R.
\end{equation}

Moreover, the total complexity is $O\left(mn\max \left\{ \sqrt{ \frac{\eig R^2}{\e \gamma}}, ~\frac{\lm_{\max}(W)m R^2}{\e^2} \right\} \right) $ arithmetic operations.
\end{theorem}
%\pd{In the Theorem 3 and below we forgot the multiplier $m$ in the complexity, see the proof of Theorem 3 in the Appendix. Please chack and update here and below.}
 
%\dd{In the experiments instead of using only one term  $p_i([\BBlm_{k+1}]_i, Y^i_1)$ for estimating $\hat{\p}_{k+1}$ on step 11 of Algorithm \ref{alg:main} we used minibatch estimation of  $\hat{\p}_{k+1}$ with size $M_k$, however the statement of the Theorem \ref{Th:WBCompl} remains correct for this case.}
We underline that even if the measures $\mu_i$, $i=1,...,m$ are discrete with large support size, it can be more efficient to apply our stochastic algorithm than to apply a deterministic algorithm. We now explain it in more details. If a measure $\mu$ is discrete, then $\W_{\gamma,\mu}^* (\blm)$ in Lemma \ref{Lm:dual_obj_properties} is represented as a finite expectation, i.e., is a sum of functions instead of an integral, and can be found explicitly. In the same way, its gradient and, hence, the gradient of the dual objective $\W^*_{\gamma}(\Blm)$ in \eqref{eq:DualObjGrad} can be found explicitly in a deterministic way. Then a deterministic accelerated primal-dual decentralized algorithm can be applied to approximate the regularized barycenter. Let us assume for simplicity that the support of measure $\mu$ is of the size $n$. Then the calculation of the exact gradient of $\W_{\gamma,\mu}^* (\blm)$ requires $O(n^2)$ arithmetic operations and the overall complexity of the deterministic algorithm is $O\left(mn^2\sqrt{\lambda_{\max}(W)R^2/\gamma\e} \right)$. For comparison, the complexity of our randomized approach in Theorem \ref{Th:WBCompl} is proportional to $n$, but not to $n^2$. So, our randomized approach is superior in the regime of large $n$. 

It is also interesting to compare the complexity of the accelerated method in Theorem \ref{Th:WBCompl} with the complexity of non-accelerated Algorithm \ref{alg:non-accler-main}. Similarly to the proof of Theorem \ref{Th:stoch_err}, extending the convergence rate proof of stochastic Mirror Descent \cite{juditsky2011solving} for the primal-dual pair of problems $(P)-(D)$, we obtain the complexity of the non-accelerated method to be $O\left(mn\max \left\{ {\lambda_{\max}(W)R^2}/{(\e\gamma)}, ~{\lm_{\max}(W) m R^2}/{\e^2} \right\} \right) $. As we see, acceleration improves the dependence on the ${\lambda_{\max}(W)R^2}/{(\e\gamma)}$, which is important, for example, for the limiting case $\gamma \to 0$, corresponding to approximation of the non-regularized barycenter.

\section{Experimental Results}\label{sec:experiments}

In this section, we present experimental results for Algorithm~\ref{alg:main}. Initially, we consider a set of agents over a network, where each agent $i$ can query realizations (i.e., samples) from a privately held random variable $Y_i \sim \mathcal{N}(\theta_i,v_i^2)$, where $\mathcal{N}(\theta,v^2)$ is a univariate Gaussian distribution with mean $\theta$  and variance $v^2$. Moreover, we set $\theta_i \in [-4,4]$ and $v_i \in [0.1,0.6]$. The objective is to compute a discrete distribution $p \in S_1(n)$ that solves \eqref{w_barycenter}. We assume $n=100$ and the support of $p$ is a set of $100$ equally spaced points on the segment $[-5,5]$. Figure~\ref{fig:performance} shows the performance of Algorithm~\ref{alg:main} for four classes of networks: complete, cycle, star, and Erd\H{o}s-R\'enyi. Moreover, we show the behavior for different network sizes, namely: $m=10,100,200,500$. Particularly we use two metrics: the function value of the dual problem and the distance to consensus, i.e., $\W^*_{\gamma}(\Blm) = \sum_{i=1}^{m} \W^*_{\gamma, \mu_i}([\BBlm]_i)$ and $C(\hat \p) := \|\sqrt{W}\hat \p\|_2$.
%\begin{align*}
%\pd{F(\BBlm) := \sum_{i=1}^{m} \W^*_{\gamma, \mu_i}([\BBlm]_i), \qquad \textrm{and} \qquad C(\hat \p) := \|\sqrt{W}\hat \p\|_2.}
%\end{align*}
As expected, when the network is a complete graph, the convergence to the final value and the distance to consensus decreases rapidly. Nevertheless, the performance in graphs with degree regularity, such as the cycle graph and the Erd\H{o}s-R\'enyi random graph, is similar to a complete graph with much less communication overhead. For the star graph, which has the worst case between the maximum and minimum number of neighbors among all nodes, the algorithms performs poorly. The reason is that despite the diameter of the graph is 2, $\lambda_{\max}(W)$, which appears in the complexity bounds, is of the order of number of vertices $m$. 

\begin{figure}[tbp!]
	\centering
	\subfigure[Local Gaussian distributions]{\begin{tikzpicture}
		\begin{axis}[
		ticklabel style = {font=\tiny},
		width=4cm,height=3.7cm,scale=0.99,
		x label style={at={(axis description cs:0.5,-0.1)},anchor=north,font=\small},
		xlabel={$x, N=1$},
%		ymode = log,
%		xmode = log,
		ymin = 0, ymax=0.06,
		xmin = -5, xmax=5, 
		cycle list name=color list,
		every axis plot/.append style={line width=1pt},
		ytick={0,0.06},
		yticklabels={$0$,$6$}],
		legend pos=south west;
		\addplot     	table [x index=40,y index=0]{ex6.dat};
		\addplot     	table [x index=40,y index=1]{ex6.dat};
		\addplot      	table [x index=40,y index=2]{ex6.dat};
		\addplot      	table [x index=40,y index=3]{ex6.dat};
		\addplot     	table [x index=40,y index=4]{ex6.dat};
		\addplot     	table [x index=40,y index=5]{ex6.dat};
		\addplot      	table [x index=40,y index=6]{ex6.dat};
		\addplot      	table [x index=40,y index=7]{ex6.dat};
		\addplot      	table [x index=40,y index=8]{ex6.dat};
		\addplot      	table [x index=40,y index=9]{ex6.dat};
		\end{axis}
		\end{tikzpicture}
\begin{tikzpicture}
	\begin{axis}[
	ticklabel style = {font=\tiny},
	width=4cm,height=3.7cm,scale=0.99,
	x label style={at={(axis description cs:0.5,-0.1)},anchor=north,font=\small},
	xlabel={$x, N=100$},
	%		ymode = log,
	%		xmode = log,
	ymin = 0, ymax=0.06,
	xmin = -5, xmax=5, 
	cycle list name=color list,
	every axis plot/.append style={line width=1pt},
	ytick={0,0.06},
	yticklabels={$0$,$6$}],
	legend pos=south west;
	\addplot     	table [x index=40,y index=10]{ex6.dat};
	\addplot     	table [x index=40,y index=11]{ex6.dat};
	\addplot      	table [x index=40,y index=12]{ex6.dat};
	\addplot      	table [x index=40,y index=13]{ex6.dat};
	\addplot     	table [x index=40,y index=14]{ex6.dat};
	\addplot     	table [x index=40,y index=15]{ex6.dat};
	\addplot      	table [x index=40,y index=16]{ex6.dat};
	\addplot      	table [x index=40,y index=17]{ex6.dat};
	\addplot      	table [x index=40,y index=18]{ex6.dat};
	\addplot      	table [x index=40,y index=19]{ex6.dat};
	\end{axis}
	\end{tikzpicture}
\begin{tikzpicture}
	\begin{axis}[
	ticklabel style = {font=\tiny},
	width=4cm,height=3.7cm,scale=0.99,
	x label style={at={(axis description cs:0.5,-0.1)},anchor=north,font=\small},
	xlabel={$x, N=200$},
	%		ymode = log,
	%		xmode = log,
	ymin = 0, ymax=0.06,
	xmin = -5, xmax=5, 
	cycle list name=color list,
	every axis plot/.append style={line width=1pt},
	ytick={0,0.06},
	yticklabels={$0$,$6$}],
	legend pos=south west;
	\addplot     	table [x index=40,y index=20]{ex6.dat};
	\addplot     	table [x index=40,y index=21]{ex6.dat};
	\addplot      	table [x index=40,y index=22]{ex6.dat};
	\addplot      	table [x index=40,y index=23]{ex6.dat};
	\addplot     	table [x index=40,y index=24]{ex6.dat};
	\addplot     	table [x index=40,y index=25]{ex6.dat};
	\addplot      	table [x index=40,y index=26]{ex6.dat};
	\addplot      	table [x index=40,y index=27]{ex6.dat};
	\addplot      	table [x index=40,y index=28]{ex6.dat};
	\addplot      	table [x index=40,y index=29]{ex6.dat};
	\end{axis}
	\end{tikzpicture}
\begin{tikzpicture}
	\begin{axis}[
	ticklabel style = {font=\tiny},
	width=4cm,height=3.7cm,scale=0.99,
	x label style={at={(axis description cs:0.5,-0.1)},anchor=north,font=\small},
	xlabel={$x, N=500$},
	%		ymode = log,
	%		xmode = log,
	ymin = 0, ymax=0.06,
	xmin = -5, xmax=5, 
	cycle list name=color list,
	every axis plot/.append style={line width=1pt},
	ytick={0,0.06},
	yticklabels={$0$,$6$}],
	legend pos=south west;
	\addplot     	table [x index=40,y index=30]{ex6.dat};
	\addplot     	table [x index=40,y index=31]{ex6.dat};
	\addplot      	table [x index=40,y index=32]{ex6.dat};
	\addplot      	table [x index=40,y index=33]{ex6.dat};
	\addplot     	table [x index=40,y index=34]{ex6.dat};
	\addplot     	table [x index=40,y index=35]{ex6.dat};
	\addplot      	table [x index=40,y index=36]{ex6.dat};
	\addplot      	table [x index=40,y index=37]{ex6.dat};
	\addplot      	table [x index=40,y index=38]{ex6.dat};
	\addplot      	table [x index=40,y index=39]{ex6.dat};
	\end{axis}
	\end{tikzpicture}
}
\\
\vspace{-0.3cm}
\subfigure[Local von Mises distributions]
{\begin{tikzpicture}
    \node at(0.2,-0.1)  {\small{$N=1$}};
		\begin{polaraxis}[
		width=3.8cm,height=3.7cm,scale=0.99,
yticklabels={,,},
    xtick = {0,90,180,270,360},
    xticklabels={0,$\pi/2$,$\pi$,$3\pi/2$},
	ticklabel style = {font=\tiny},
		cycle list name=color list,
		data cs = polarrad,
		every axis plot/.append style={line width=1pt}
		]
		\addplot      	table [x index=40,y index=0]{ex7.dat};
		\addplot       	table [x index=40,y index=1]{ex7.dat};
		\addplot      	table [x index=40,y index=2]{ex7.dat};
		\addplot       	table [x index=40,y index=3]{ex7.dat};
		\addplot      	table [x index=40,y index=4]{ex7.dat};
		\addplot       	table [x index=40,y index=5]{ex7.dat};
		\addplot      	table [x index=40,y index=6]{ex7.dat};
		\addplot       	table [x index=40,y index=7]{ex7.dat};
		\addplot      	table [x index=40,y index=8]{ex7.dat};
		\addplot       	table [x index=40,y index=9]{ex7.dat};
		\end{polaraxis}
		\end{tikzpicture}
\begin{tikzpicture}
    \node at (0.2,-0.1)  {\small $N=100$};
	\begin{polaraxis}[
	width=3.8cm,height=3.7cm,scale=0.99,
	yticklabels={,,},
    xtick = {0,90,180,270,360},
    xticklabels={0,$\pi/2$,$\pi$,$3\pi/2$},
	ticklabel style = {font=\tiny},
	cycle list name=color list,
	data cs = polarrad,
	every axis plot/.append style={line width=1pt}
	]
	\addplot      	table [x index=40,y index=10]{ex7.dat};
	\addplot       	table [x index=40,y index=11]{ex7.dat};
	\addplot      	table [x index=40,y index=12]{ex7.dat};
	\addplot       	table [x index=40,y index=13]{ex7.dat};
	\addplot      	table [x index=40,y index=14]{ex7.dat};
	\addplot       	table [x index=40,y index=15]{ex7.dat};
	\addplot      	table [x index=40,y index=16]{ex7.dat};
	\addplot       	table [x index=40,y index=17]{ex7.dat};
	\addplot      	table [x index=40,y index=18]{ex7.dat};
	\addplot       	table [x index=40,y index=19]{ex7.dat};
	\end{polaraxis}
	\end{tikzpicture}
\begin{tikzpicture}
    \node at (0.2,-0.1)  {\small $N=200$};
	\begin{polaraxis}[
	width=3.8cm,height=3.7cm,scale=0.99,
	yticklabels={,,},
    xtick = {0,90,180,270,360},
    xticklabels={0,$\pi/2$,$\pi$,$3\pi/2$},
	ticklabel style = {font=\tiny},
	cycle list name=color list,
	data cs = polarrad,
	every axis plot/.append style={line width=1pt}
	]
	\addplot      	table [x index=40,y index=20]{ex7.dat};
	\addplot       	table [x index=40,y index=21]{ex7.dat};
	\addplot      	table [x index=40,y index=22]{ex7.dat};
	\addplot       	table [x index=40,y index=23]{ex7.dat};
	\addplot      	table [x index=40,y index=24]{ex7.dat};
	\addplot       	table [x index=40,y index=25]{ex7.dat};
	\addplot      	table [x index=40,y index=26]{ex7.dat};
	\addplot       	table [x index=40,y index=27]{ex7.dat};
	\addplot      	table [x index=40,y index=28]{ex7.dat};
	\addplot       	table [x index=40,y index=29]{ex7.dat};
	\end{polaraxis}
	\end{tikzpicture}
\begin{tikzpicture}
    \node at (0.2,-0.1)  {\small $N=500$};
	\begin{polaraxis}[
	width=3.8cm,height=3.7cm,scale=0.99,
	yticklabels={,,},
	%xticklabel={
	%	\pgfmathparse{\tick/180}
	%	\pgfmathifisint{\pgfmathresult}{$\pgfmathprintnumber[int detect]{\pgfmathresult}\pi$}%
	%	{$\pgfmathprintnumber[frac,frac denom=6,frac whole=false]{\pgfmathresult}\pi$}
	%},
    xtick = {0,90,180,270,360},
    xticklabels={0,$\pi/2$,$\pi$,$3\pi/2$},
	ticklabel style = {font=\tiny},
	cycle list name=color list,
	data cs = polarrad,
	every axis plot/.append style={line width=1pt}
	]
	\addplot      	table [x index=40,y index=30]{ex7.dat};
	\addplot       	table [x index=40,y index=31]{ex7.dat};
	\addplot      	table [x index=40,y index=32]{ex7.dat};
	\addplot       	table [x index=40,y index=33]{ex7.dat};
	\addplot      	table [x index=40,y index=34]{ex7.dat};
	\addplot       	table [x index=40,y index=35]{ex7.dat};
	\addplot      	table [x index=40,y index=36]{ex7.dat};
	\addplot       	table [x index=40,y index=37]{ex7.dat};
	\addplot      	table [x index=40,y index=38]{ex7.dat};
	\addplot       	table [x index=40,y index=39]{ex7.dat};
	\end{polaraxis}
	\end{tikzpicture}
}
\vspace{-0.3cm}
	\caption{\small Local barycenter of a set of Gaussian distribution and von Mises distributions. Barycenter is generated by the Algorithm~\ref{alg:main} for a set of $10$ agents over an Erd\H{o}s-R\'enyi random graph at different iteration numbers. %(a) Each agent can access private realizations from a Gaussian random variable. (b) 
    Each agent can access private realizations from a von Mises random variable.}
	\label{fig:gaussians}
\end{figure}
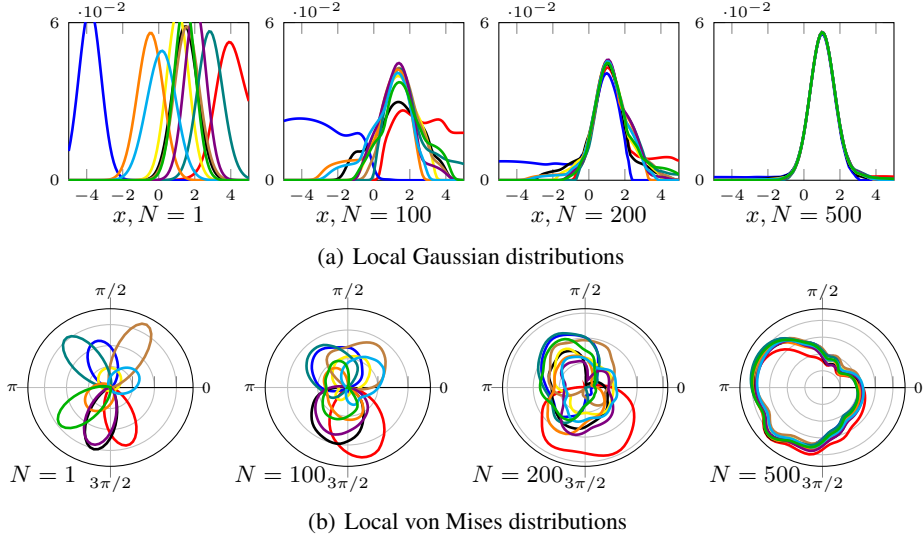

Figure~\ref{fig:gaussians}(a) shows a sample of the local barycenters of $10$ agents on an Erd\H{o}s-R\'enyi random graph, with local Gaussian distributions, at different times of the Algorithm~\ref{alg:main}, $N=1,100,200,500$. The local barycenters of all the agents in the network converge to a common distribution. Similarly, Figure~\ref{fig:gaussians}(b) shows the convergence of the local barycenters of the agents on the same Erd\H{o}s-R\'enyi random graph when the local distributions are von Mises distributions. Particularly, for the cases of von Mises distributions, we have used the angle between to points distance function.

\begin{figure}[tbp!]
	\centering
	\resizebox{10cm}{!}{
		\begin{tikzpicture}
		\draw (0,0) -- (11,0) -- (11,0.5) -- (0,0.5) -- (0,0);
		\draw[red,line width=2pt] (0.1,0.25) -- (0.6,0.25) node[black] at (1.3,0.25) {$\text{Cycle}$};
		\draw[blue,line width=2pt] (3,0.25) -- (3.5,0.25) node[black] at (4.6,0.25) {$\text{Erd\H{o}s-R\'enyi}$};
		\draw[green,line width=2pt] (5.9,0.25) -- (6.4,0.25) node[black] at (7.1,0.25) {$\text{Star}$};
		\draw[black,line width=2pt] (8.8,0.25) -- (9.3,0.25) node[black] at (10.0,0.25) {$\text{Complete}$};
		\end{tikzpicture}
	}\\
    \vspace{-0.2cm}
	\subfigure{\begin{tikzpicture}
		\begin{axis}[
		ticklabel style = {font=\tiny},
		width=3.8cm,height=3.5cm,scale=0.99,
		x label style={at={(axis description cs:0.5,-0.1)},anchor=north,font=\small},
		y label style={at={(axis description cs:-0.1,.5)},anchor=south},
		xlabel={Iterations $m=200$},
		ylabel={$F(\tilde \lambda_k)$},
		%		ymode = log,
		%		xmode = log,
		ymin = -2, ymax=1,
		xmin = 1, xmax=1000, 
		every axis plot/.append style={line width=1pt}],
		legend pos=south west;
		\addplot [black]     	table [x index=8,y index=0]{ex1.dat};
		\addplot [red]     		table [x index=8,y index=2]{ex1.dat};
		\addplot [green]     	table [x index=8,y index=4]{ex1.dat};
		\addplot [blue]     	table [x index=8,y index=6]{ex1.dat};
		\end{axis}
		\end{tikzpicture}
	}
	\subfigure{\begin{tikzpicture}
		\begin{axis}[
		ticklabel style = {font=\tiny},
		width=3.8cm,height=3.5cm,scale=0.99,
		x label style={at={(axis description cs:0.5,-0.1)},anchor=north,font=\small},
		xlabel={Iterations $m=100$},
		%	ylabel={$F(\tilde \lambda_k)$},
		%		ymode = log,
		%		xmode = log,
		ymin = -2, ymax=1,
		xmin = 1, xmax=1000, 
		every axis plot/.append style={line width=1pt}],
		legend pos=south west;
		\addplot [black]     	table [x index=8,y index=0]{ex2.dat};
		\addplot [red]     		table [x index=8,y index=2]{ex2.dat};
		\addplot [green]     	table [x index=8,y index=4]{ex2.dat};
		\addplot [blue]     	table [x index=8,y index=6]{ex2.dat};
		\end{axis}
		\end{tikzpicture}
	}
	\subfigure{\begin{tikzpicture}
		\begin{axis}[
		ticklabel style = {font=\tiny},
		width=3.8cm,height=3.5cm,scale=0.99,
		x label style={at={(axis description cs:0.5,-0.1)},anchor=north,font=\small},
		xlabel={Iterations $m=10$},
		%	ylabel={$F(\tilde \lambda_k)$},
		%		ymode = log,
		%		xmode = log,
		ymin = -2, ymax=1,
		xmin = 1, xmax=1000, 
		every axis plot/.append style={line width=1pt}],
		legend pos=south west;
		\addplot [black]     	table [x index=8,y index=0]{ex3.dat};
		\addplot [red]     		table [x index=8,y index=2]{ex3.dat};
		\addplot [green]     	table [x index=8,y index=4]{ex3.dat};
		\addplot [blue]     	table [x index=8,y index=6]{ex3.dat};
		\end{axis}
		\end{tikzpicture}
	}
	\subfigure{\begin{tikzpicture}
		\begin{axis}[
		ticklabel style = {font=\tiny},
		width=3.8cm,height=3.5cm,scale=0.99,
		x label style={at={(axis description cs:0.5,-0.1)},anchor=north,font=\small},
		xlabel={Iterations $m=500$},
		%	ylabel={$F(\tilde \lambda_k)$},
		%		ymode = log,
		%		xmode = log,
		ymin = -2, ymax=1,
		xmin = 1, xmax=1000, 
		every axis plot/.append style={line width=1pt}],
		legend pos=south west;
		\addplot [black]     	table [x index=8,y index=0]{ex4.dat};
		\addplot [red]     		table [x index=8,y index=2]{ex4.dat};
		\addplot [green]     	table [x index=8,y index=4]{ex4.dat};
		\addplot [blue]     	table [x index=8,y index=6]{ex4.dat};
		\end{axis}
		\end{tikzpicture}
	}
	\\
    \vspace{-0.4cm}
	\subfigure{\begin{tikzpicture}
		\begin{axis}[
		ticklabel style = {font=\tiny},
		width=3.8cm,height=3.5cm,scale=0.99,
		x label style={at={(axis description cs:0.5,-0.1)},anchor=north,font=\small},
		xlabel={Iterations $m=200$},
		ylabel={$C(\hat p_k)$},
		%		ymode = log,
		%		xmode = log,
		ymin = 0, ymax=1,
		xmin = 1, xmax=1000, 
		every axis plot/.append style={line width=1pt}],
		legend pos=south west;
		\addplot [black]     	table [x index=8,y index=1]{ex4.dat};
		\addplot [red]     		table [x index=8,y index=3]{ex4.dat};
		\addplot [green]     	table [x index=8,y index=5]{ex4.dat};
		\addplot [blue]     	table [x index=8,y index=7]{ex4.dat};
		\end{axis}
		\end{tikzpicture}
	}
	\subfigure{\begin{tikzpicture}
		\begin{axis}[
		ticklabel style = {font=\tiny},
		width=3.8cm,height=3.5cm,scale=0.99,
		x label style={at={(axis description cs:0.5,-0.1)},anchor=north,font=\small},
		xlabel={Iterations $m=100$},
		%	ylabel={$C(\hat p_k)$},
		%		ymode = log,
		%		xmode = log,
		ymin = 0, ymax=1,
		xmin = 1, xmax=1000, 
		every axis plot/.append style={line width=1pt}],
		legend pos=south west;
		\addplot [black]     	table [x index=8,y index=1]{ex1.dat};
		\addplot [red]     		table [x index=8,y index=3]{ex1.dat};
		\addplot [green]     	table [x index=8,y index=5]{ex1.dat};
		\addplot [blue]     	table [x index=8,y index=7]{ex1.dat};
		\end{axis}
		\end{tikzpicture}
	}
	\subfigure{\begin{tikzpicture}
		\begin{axis}[
		ticklabel style = {font=\tiny},
		width=3.8cm,height=3.5cm,scale=0.99,
		x label style={at={(axis description cs:0.5,-0.1)},anchor=north,font=\small},
		xlabel={Iterations $m=10$},
		%	ylabel={$C(\hat p_k)$},
		%		ymode = log,
		%		xmode = log,
		ymin = 0, ymax=1,
		xmin = 1, xmax=1000, 
		every axis plot/.append style={line width=1pt}],
		legend pos=south west;
		\addplot [black]     	table [x index=8,y index=1]{ex2.dat};
		\addplot [red]     		table [x index=8,y index=3]{ex2.dat};
		\addplot [green]     	table [x index=8,y index=5]{ex2.dat};
		\addplot [blue]     	table [x index=8,y index=7]{ex2.dat};
		\end{axis}
		\end{tikzpicture}
	}
	\subfigure{\begin{tikzpicture}
		\begin{axis}[
		ticklabel style = {font=\tiny},
		width=3.8cm,height=3.5cm,scale=0.99,
		x label style={at={(axis description cs:0.5,-0.1)},anchor=north,font=\small},
		xlabel={Iterations $m=500$},
		%	ylabel={$C(\hat p_k)$},
		%		ymode = log,
		%		xmode = log,
		ymin = 0, ymax=1,
		xmin = 1, xmax=1000, 
		every axis plot/.append style={line width=1pt}],
		legend pos=south west;
		\addplot [black]     	table [x index=8,y index=1]{ex3.dat};
		\addplot [red]     		table [x index=8,y index=3]{ex3.dat};
		\addplot [green]     	table [x index=8,y index=5]{ex3.dat};
		\addplot [blue]     	table [x index=8,y index=7]{ex3.dat};
		\end{axis}
		\end{tikzpicture}
	}
    \vspace{-0.6cm}
	\caption{\small Dual function value and distance to consensus for $200,100,10,500$ agents, $M_k = 100$ and $\gamma =  0.1$.}
	\label{fig:performance}
\end{figure}
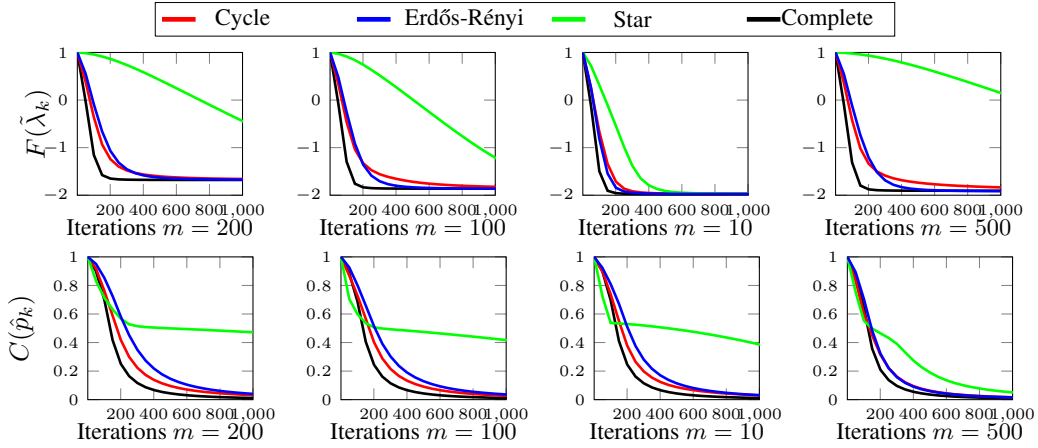

\begin{figure}[tb!]
    \centering
    \fbox{\includegraphics[width=0.17\textwidth]{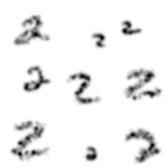}}
    \fbox{\includegraphics[width=0.17\textwidth]{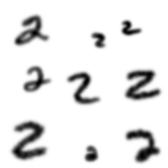}}
    \fbox{\includegraphics[width=0.17\textwidth]{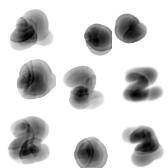}}
    \fbox{\includegraphics[width=0.17\textwidth]{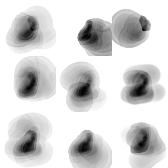}}
    \fbox{\includegraphics[width=0.17\textwidth]{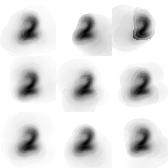}}\\
        \begin{tikzpicture}
    \draw node[black] at (0.3,0.25) {$N=1$};
    \draw node[black] at (2.7,0.25) {$N=1000$};
    \draw node[black] at (5.4,0.25) {$N=2000$};
    \draw node[black] at (8.0,0.25) {$N=3000$};
    \draw node[black] at (10.5,0.25) {$N=4000$};
    \end{tikzpicture}
    \vspace{-0.4cm}
    \caption{\small Wasserstein barycenter of a subset of images of the digit $2$ from the MNIST dataset~\cite{LeCun1998}. Each block shows a subset of $9$ randomly selected local barycenters, generated by Algorithm~\ref{alg:main} at different time instances. The $9$ agents are selected from a network of $500$ agents on an Erd\H{o}s-R\'enyi random graph.}
    \label{fig:numbers}
\end{figure}
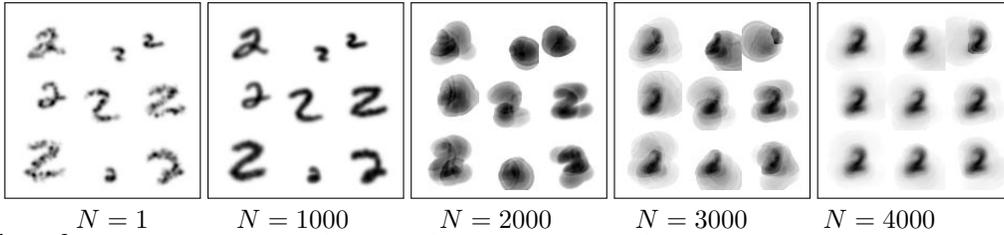

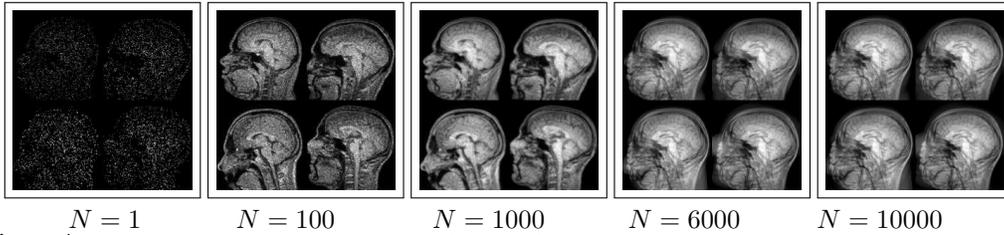
\begin{figure}[tb!]
    \centering
    \fbox{\includegraphics[width=0.17\textwidth]{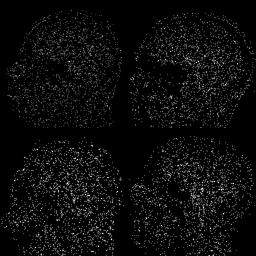}}
    \fbox{\includegraphics[width=0.17\textwidth]{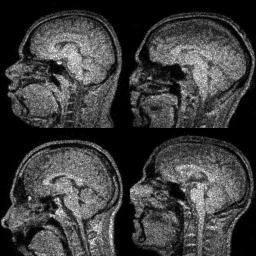}}
    \fbox{\includegraphics[width=0.17\textwidth]{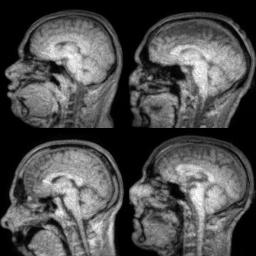}}
    \fbox{\includegraphics[width=0.17\textwidth]{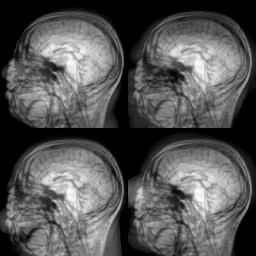}}
    \fbox{\includegraphics[width=0.17\textwidth]{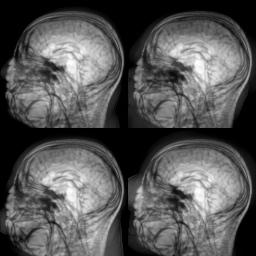}}\\
    \begin{tikzpicture}
    \draw node[black] at (0.3,0.25) {$N=1$};
    \draw node[black] at (2.7,0.25) {$N=100$};
    \draw node[black] at (5.4,0.25) {$N=1000$};
    \draw node[black] at (8.0,0.25) {$N=6000$};
    \draw node[black] at (10.6,0.25) {$N=10000$};
    \end{tikzpicture}
    \vspace{-0.4cm}
\caption{\small Wasserstein barycenter for a subset of images from the IXI dataset~\cite{ixidata}. Each block shows the local barycenters, of $4$ agents, generated by Algorithm~\ref{alg:main} at different time instances. The $4$ agents are connected on a cycle graph.}
    \label{fig:MRI}
\end{figure}

Figure~\ref{fig:numbers} shows the computed local barycenter of $9$ agents in a network of $500$ nodes at different iteration numbers. Each agent holds a local copy of a sample of the digit $2$ ($56 \times 56$ image) from the MNIST dataset~\cite{LeCun1998}. All agents converge to the same image that structurally represents the aggregation of the original $500$ images held over the network. 
Finally, Figure~\ref{fig:MRI} shows a simple example of an application of Wasserstein barycenter on medical image aggregation where we have $4$ agents connected over a cycle graph and each agent holds a magnetic resonance image ($256 \times 256$) from the IXI dataset~\cite{ixidata}.

\section{Conclusions and Future Directions}\label{sec:conclusions}

We propose a novel distributed algorithm for the computation of the regularized Wasserstein barycenter of a set of continuous measures stored distributedly over a network of agents. Moreover, we provide explicit and non-asymptotic iteration and sample complexity analysis in terms of the problem parameters and the network topology. Our algorithm is based on a new general algorithm for the solution of stochastic convex optimization problems with linear constraints. In contrast to the recent literature, our algorithm can be executed over arbitrary connected and static networks where nodes are oblivious to the network topology, which makes it suitable for large-scale network optimization setting. Additionally, our analysis indicates that the randomization strategy provides faster convergence rates than the deterministic procedure when the support size of the barycenter is large.

The presented experiments were carried out in a single machine and %Sometimes, off-the-shelf parallel toolboxes were used, but not with the intent of recreating a real large-scale distributed scenario.
implementation of our algorithm on real networks is a major research thrust for future projects. 
Extending fast distributed algorithms for the case of time-varying and directed graph networks remains an open question. Notably, it is not clear what is the effect of the network dynamics in the quality of the solution of specific problems such as the Wasserstein barycenter. Moreover, efficient communication strategies between nodes should be considered as well. 
The extension to the decentralized distributed setting of Sinkhorn-type algorithms \cite{benamou2015iterative} for regularized Wasserstein barycenter and other related algorithms, e.g., Wasserstein propagation
\cite{solomon2015convolutional}, requires further work.

\subsubsection*{Acknowledgments}
The work of A. Nedi\'c and C.A. Uribe in Sect. 5 is supported by the National Science Foundation under grant
no. CPS 15-44953. The research by P. Dvurechensky, D. Dvinskikh, and A. Gasnikov in Sect. 3 and Sect. 4 was funded by the Russian Science Foundation (project 18-71-10108). 
\newpage

\medskip
\small

\bibliographystyle{abbrv}
%\bibliographystyle{achemso}

%\bibliography{wass,PD_references,time_varying}

% \setcounter{Lm}{0}
% \renewcommand{\theLm}{\thesection.\arabic{Lm}}
% \renewcommand{\theTh}{\thesection.\arabic{Th}}

\newpage
\appendix
\section{Proofs and Additional Numerical Results}

In this appendix, we present the complete proofs of the Lemmas and Theorems stated in the main article. Moreover, we show additional experimental results. The contents of the Appendix are organized as follows:

\begin{itemize}
\item Subsection \ref{supp:lemma_1}, Subsection \ref{supp:the_1}, Subsection \ref{supp:the_2}, Subsection \ref{supp:lemma_2}, and Subsection~\ref{supp:the_3} present the complete proofs of the Lemmas and Theorems of the main paper.
\item Subsection \ref{supp:networks} shows a graphic representation of the network topologies used for the experimental results, namely: complete graph, star graph, cycle graph and Erd\H{o}s-R\'enyi random graph.
\item Subsection \ref{supp:gauss} shows, for various time instances, the local Wasserstein barycenter of $10$ agents connected on an Erd\H{o}s-R\'enyi random graph. Each agent holds a private Gaussian measure from which it can query samples. Different colors represent different agents. Time evolves with the number of iterations.
\item Subsection \ref{supp:von} shows, for various time instances, local Wasserstein barycenter of $10$ agents connected on an Erd\H{o}s-R\'enyi random graph. Each agent holds a private von Mises measure from which it can query samples. Different colors represent different agents. Time evolves with the number of iterations. 
\item Subsection \ref{supp:mnist}: shows, for various time instances, local Wasserstein barycenter of $100$ agents connected on an Erd\H{o}s-R\'enyi random graph. Each agent holds a private sample of the digit $2$ from the MNIST dataset. We assume the normalize image as a probability distribution from which agents can sample from. Time evolves with the number of iterations.
\item Subsection \ref{supp:ixi}: shows, for various time instances, local Wasserstein barycenter of $4$ agents connected on an cycle graph. Each agent holds a private sample of an magnetic resonance image from the IXI dataset. We assume the normalize image as a probability distribution from which agents can sample from. Time evolves with the number of iterations.
\item Attached videos
\begin{itemize}
\item \texttt{Gauss\_ex1.avi}: Example $1$. The local Wasserstein barycenter of $10$ agents connected on an Erd\H{o}s-R\'enyi random graph. Each agent holds a private Gaussian measure from which it can query samples. Different colors represent different agents. Time evolves with the number of iterations. 
\item \texttt{Gauss\_ex2.avi}: Example $2$. The local Wasserstein barycenter of $10$ agents connected on an Erd\H{o}s-R\'enyi random graph. Each agent holds a private Gaussian measure from which it can query samples. Different colors represent different agents. Time evolves with the number of iterations. 
\item \texttt{MNIST\_digit2.avi}: The local Wasserstein barycenter of $100$ agents connected on an Erd\H{o}s-R\'enyi random graph. Each agent holds a private sample of the digit $2$ from the MNIST dataset. We assume the normalize image as a probability distribution from which agents can sample from. Time evolves with the number of iterations.
\item \texttt{MNIST\_digit3.avi}: The local Wasserstein barycenter of $100$ agents connected on an Erd\H{o}s-R\'enyi random graph. Each agent holds a private sample of the digit $3$ from the MNIST dataset. We assume the normalize image as a probability distribution from which agents can sample from. Time evolves with the number of iterations.
\item \texttt{von\_mises\_ex1.avi}: Example $1$. The local Wasserstein barycenter of $10$ agents connected on an Erd\H{o}s-R\'enyi random graph. Each agent holds a private von Mises measure from which it can query samples. Different colors represent different agents. Time evolves with the number of iterations. 
\item \texttt{von\_mises\_ex2.avi}: Example $2$. The local Wasserstein barycenter of $10$ agents connected on an Erd\H{o}s-R\'enyi random graph. Each agent holds a private von Mises measure from which it can query samples. Different colors represent different agents. Time evolves with the number of iterations. 
\item \texttt{ixi\_mr.avi}: The local Wasserstein barycenter of $4$ agents connected on an cycle graph. Each agent holds a private sample of an magnetic resonance image from the IXI dataset. We assume the normalize image as a probability distribution from which agents can sample from. Time evolves with the number of iterations.
\end{itemize}
\end{itemize}

\subsection{Proof of Lemma \ref{Lm:dual_obj_properties}}\label{supp:lemma_1}

%For entropy-regularized dual formulation of OT \eqref{dual_OT} we can explicitly write the notion of c-transform.
%\begin{align}\label{c_trans} 
%\lambda^{c,\gamma}(y) = -\gamma\log\left(\sum_{i=1}^n \exp\left(\frac{\lambda_i -c(y,z_i)}{\gamma}\right)p_i \right)
%\end{align}
%Using this c-transform \eqref{dual_OT} is
%equivalent to the following optimization problem
%\begin{align}\label{conj}
%\mathcal{W}_{\gamma,\mu}(p) &=\max_\lambda \left\{\int_\mathcal{Y}\lambda^{c,\gamma}(y)q(y)dy + \sum_{i=1}^n\lambda_ip_i\right\} = \max_\lambda \left\{\int_\mathcal{Y}\lambda^{c,\gamma}(y)q(y)dy + \la \lambda,p\ra \right\} \notag \\
%&= \max_\lambda \left\{-\mathcal{W}^*_\gamma (\lambda)+ \la \lambda,p\ra \right\}
%\end{align}
Primal and dual optimal transport problem corresponding to the regularized Wasserstein distance can be written as follows
{\small \begin{align*}
\dd{\W_{\gamma,\mu}(p)}&=\min_{\pi \in \Pi(\mu,\nu)}\left\{\sum_{l=1}^n\int_{\Y}  c_l(y)\pi_l(y)dy +\gamma \sum_{l=1}^n\int_{\Y} \pi_l(y)\log \pi_l(y)dy - \gamma\log\xi \right\}\\
&= \max_{\dd{\bar{\lambda}} \in \R^n,v \in \C(\X)}\left\{\sum_{l=1}^n [p]_l [\bar{\lambda}]_l +\int_{\Y}q(y)v(y)dy -\gamma\sum_{l=1}^n\int_{\Y}\exp\left(\frac{\dd{[\bar{\lambda}}]_l-c_l(y)+v(y)}{\gamma} -1\right)dy  \right\} \\
&= \max_{\dd{\bar{\lambda}}\in \R^n} \left\{\la p, \dd{\bar{\lambda}}\ra - \gamma\int_\Y \log \left( \frac{1}{q(y)}\sum_{l=1}^n\exp\left(\frac{\dd{[\bar{\lambda}]}_l-c_l(y)}{\gamma}\right) \right)q(y)dy \right\},
\end{align*}}
where we used that $\xi$ is the uniform distribution on $\Y \times \mathcal{Z}$. By the definition of Fenchel-Legendre transform, using that $\W_{\gamma,\mu}(p) = (\W_{\gamma,\mu}(p))^{**}$, we get the first statement of the Lemma
{\small \begin{align*}
\W^*_{\gamma,\mu} (\bar{\lambda}) 
&= \gamma \int_{\mathcal{Y}} \log\left(\frac{1}{q(y)}\sum_{l=1}^n \exp\left( ([\bar{\lambda}]_l -c_l(y))/{\gamma}\right)\right)q(y)dy \\
\\
&=\E_{Y\sim\mu}\gamma \log\left(\frac{1}{q(Y)}\sum_{l=1}^n\exp\left( ([\blm]_l-c_l(Y))/{\gamma} \right)\right),
\end{align*}}
where $Y\sim \mu$ means that random variable $Y$ distributed according to measure $\mu$.

Differentiating, we obtain that the $l$-th component of the gradient of $\W^*_{\gamma,\mu} (\bar{\lambda})$ is
 {\small \begin{align*}
[\nabla \W_{\gamma,\mu}^*(\bar{\lambda})]_l 
&= \int_{\mathcal{Y}} \frac{\exp(([\bar{\lambda}]_l-c_l(y))/\gamma)  }{\sum_{\ell=1}^n\exp(([\bar{\lambda}]_\ell-c_\ell(y))/\gamma)} q(y)dy\\
\\
&=\E_{Y\sim\mu} \frac{\exp(([\bar{\lambda}]_l-c_l(Y))/\gamma)  }{\sum_{\ell=1}^n\exp(([\bar{\lambda}]_\ell-c_\ell(Y))/\gamma)} , ~l=1,\dots,n.
\end{align*}}
\noindent To prove the Lipschitz continuity of this gradient, we calculate the diagonal elements of the Hessian 
{\small \begin{align*}
[\nabla^2 \W^*_{\gamma,\mu}(\bar{\lambda})]_{ll} 
&= \frac{1}{\gamma}\int_\Y \frac{\exp(([\bar{\lambda}]_l-c_l(y))/\gamma)\sum_{\ell=1}^n\exp(([\bar{\lambda}]_\ell-c_\ell(y))/\gamma) - \exp^2(([\bar{\lambda}]_l-c_l(y))/\gamma)  }{\left(\sum_{\ell=1}^n\exp(([\bar{\lambda}]_\ell-c_\ell(y))/\gamma)\right)^2}q(y)dy
\end{align*}}
and estimate its trace
\begin{align*}
{\rm Tr}(\nabla^2 \W^*_{\gamma,\mu}(\bar{\lambda})) &\leq \frac{1}{\gamma}\int_\Y \frac{\sum_{l=1}^n\exp(([\bar{\lambda}]_l-c_l(y))/\gamma)\sum_{\ell=1}^n\exp(([\bar{\lambda}]_\ell-c_\ell(y))/\gamma)}{\left(\sum_{\ell=1}^n\exp(([\bar{\lambda}]_\ell-c_\ell(y))/\gamma)\right)^2}q(y)dy \\
&= \frac{1}{\gamma}\int_\Y q(y)dy = \frac{1}{\gamma}. 
\end{align*}
This inequality proves that $\nabla \W_{\gamma,\mu}^*(\bar{\lambda})$ is $\frac{1}{\gamma}$- Lipschitz continuous with respect to the 2-norm.

\subsection{Proof of Theorem~\ref{Th:ASGDConv}}\label{supp:the_1}

Let us fix an arbitrary $\lambda \in \Lambda$. From the optimality condition in \eqref{eq:zetakp1Def}, we have
\begin{equation}
\la \nabla V[\zeta_{k}](\zeta_{k+1}) + \alpha_{k+1} \nabla\Phi(\lambda_{k+1},\xi_{k+1}), \lambda - \zeta_{k+1} \ra \geq 0.
\label{eq:Lm1Pr1}
\end{equation}
Further, 
\begin{align*}
\alpha_{k+1}\la \nabla\Phi(\lambda_{k+1},\xi_{k+1}), \zeta_{k} - \lambda\ra &= \\
 &\hspace{-11em}= \alpha_{k+1}\la \nabla \Phi(\lambda_{k+1}, \xi_{k+1}), \zeta_{k} - \zeta_{k+1}\ra + \alpha_{k+1}\la \nabla\Phi(\lambda_{k+1},\xi_{k+1}), \zeta_{k+1} - \lambda\ra \notag \\
    &\hspace{-11em}\stackrel{\eqref{eq:Lm1Pr1}}{\leq} \alpha_{k+1}\la  \dd{\nabla \Phi(\lambda_{k+1},\xi_{k+1})}, \zeta_{k} - \zeta_{k+1}\ra + \la - \nabla V[\zeta_{k}](\zeta_{k+1}) , \zeta_{k+1} - \lambda  \ra \notag \\
	  &\hspace{-11em}= \alpha_{k+1}\la \nabla\Phi(\lambda_{k+1},\xi_{k+1}), \zeta_{k} - \zeta_{k+1}\ra + V[\zeta_k](\lambda) - V[\zeta_{k+1}](\lambda) - V[\zeta_k](\zeta_{k+1}) \notag \\
		&\hspace{-11em}\stackrel{\eqref{eq:BFLowBound}}{\leq} \alpha_{k+1}\la \nabla\Phi(\lambda_{k+1},\xi_{k+1}), \zeta_{k} - \zeta_{k+1}\ra + V[\zeta_k](\lambda) - V[\zeta_{k+1}](\lambda) - \frac12\|\zeta_k-\zeta_{k+1}\|^2 \notag \\
&\hspace{-13em}\stackrel{\eqref{eq:lambdakp1Def},\eqref{eq:etakp1Def}}{=} C_{k+1} \la \nabla\Phi(\lambda_{k+1},\xi_{k+1}), \lambda_{k+1} - \eta_{k+1}\ra + V[\zeta_k](\lambda) - V[\zeta_{k+1}](\lambda)- \frac{C_{k+1}^2}{2\alpha_{k+1}^2}\|\lambda_{k+1} - \eta_{k+1}\|^2 \notag \\	&\hspace{-11em}\stackrel{\eqref{eq:alpQuadEq}}{=} C_{k+1}\left(\la \nabla\Phi(\lambda_{k+1},\xi_{k+1}), \lambda_{k+1} - \eta_{k+1}\ra - \frac{2L}{2}\|\lambda_{k+1} - \eta_{k+1}\|^2  \right) + V[\zeta_k](\lambda) - V[\zeta_{k+1}](\lambda). \notag
\end{align*}
Add and subtract the term $C_{k+1}\la \nabla \vp(\lambda_{k+1}), \lambda_{k+1} - \eta_{k+1}\ra$, then
\begin{align}\label{eq1}
\alpha_{k+1}\la\nabla\Phi(\lambda_{k+1},\xi_{k+1}), \zeta_{k} - \lambda\ra 
&\leq C_{k+1}\left(\la \nabla\Phi(\lambda_{k+1},\xi_{k+1})- \nabla \vp(\lambda_{k+1}), \lambda_{k+1} - \eta_{k+1}\ra \right.\notag\\
&\hspace{-7em} \left. - \frac{2L}{2}\|\lambda_{k+1} - \eta_{k+1}\|^2 + \la \nabla \vp(\lambda_{k+1}), \lambda_{k+1} - \eta_{k+1}\ra \right) + V[\zeta_k](\lambda) - V[\zeta_{k+1}](\lambda). 
\end{align}
Using Fenchel inequality $\langle g,x \rangle \leq \frac{1}{2\zeta}\|g\|_{*}^2 + \frac{\zeta}{2}\|x\|^2$, we estimate 
\begin{align*}
\la \nabla\Phi(\lambda_{k+1},\xi_{k+1})- \nabla \vp(\lambda_{k+1}), \lambda_{k+1} - \eta_{k+1}\ra &\leq\\
&\hspace{-10em} \leq \frac{1}{2L}\| \nabla\Phi(\lambda_{k+1},\xi_{k+1})- \nabla \vp(\lambda_{k+1})\|_{*}^2 + \frac{L}{2}\| \lambda_{k+1} - \eta_{k+1}\|^2.
\end{align*}
Therefore, we can rewrite \eqref{eq1} as
\begin{align}\label{eq_lemma1}
\alpha_{k+1}\la \nabla\Phi(\lambda_{k+1},\xi_{k+1}), \zeta_{k} - \lambda\ra &\leq \notag\\
&\hspace{-11em}\leq C_{k+1}\left(\frac{1}{2L}\| \nabla\Phi(\lambda_{k+1},\xi_{k+1})- \nabla \vp(\lambda_{k+1})\|_{*}^2 \right.\notag\\
&\hspace{-9em}+ \left.\la \nabla \vp(\lambda_{k+1}), \lambda_{k+1} - \eta_{k+1}\ra - \frac{L}{2}\|\lambda_{k+1} - \eta_{k+1}\|^2  \right) + V[\zeta_k](\lambda) - V[\zeta_{k+1}](\lambda) \notag\\
&\hspace{-11em} = C_{k+1}\left(\la \nabla \vp(\lambda_{k+1}), \lambda_{k+1} - \eta_{k+1}\ra - \frac{L}{2}\|\lambda_{k+1} - \eta_{k+1}\|^2  \right) + V[\zeta_k](\lambda) - V[\zeta_{k+1}](\lambda) \notag\\
&\hspace{-9em} +  \frac{C_{k+1}}{2L}\| \nabla\Phi(\lambda_{k+1},\xi_{k+1})- \nabla \vp(\lambda_{k+1})\|_{*}^2 \notag\\
&\hspace{-11em}\stackrel{\eqref{eq:nfLipDef}}{\leq} C_{k+1}\left(  \vp(\lambda_{k+1}) - \vp(\eta_{k+1})  \right)  + V[\zeta_k](\lambda) - V[\zeta_{k+1}](\lambda) \notag\\ 
&\hspace{-9em}+\frac{C_{k+1}}{2L}\| \nabla \Phi(\lambda_{k+1},\xi_{k+1})- \nabla \vp(\lambda_{k+1})\|_{*}^2.  
\end{align}
Similarly, adding and substracting the term $\la \nabla \vp(\lambda_{k+1}), \eta_{k} - \lambda_{k+1}\ra$, we have
\begin{align}\label{eq_lemma2}
\hspace{-2em}\la \nabla\Phi(\lambda_{k+1},\xi_{k+1}), \eta_{k} - \lambda_{k+1} \ra 
&\leq \notag \\
&\hspace{-6em} \leq \la \nabla \vp(\lambda_{k+1}), \eta_{k} - \lambda_{k+1} \ra + \la \nabla\Phi(\lambda_{k+1},\xi_{k+1}) - \nabla \vp(\lambda_{k+1}), \eta_{k} - \lambda_{k+1} \ra  \notag\\
&\hspace{-6em}\stackrel{\text{conv-ty}}{\leq}  \vp(\eta_k) - \vp(\lambda_{k+1}) + \la \nabla \Phi(\lambda_{k+1}, \xi_{k+1}) - \nabla\vp(\lambda_{k+1}), \eta_{k} - \lambda_{k+1} \ra.
\end{align}
Finally, 
\begin{align*}
\alpha_{k+1}\la \nabla \Phi(\lambda_{k+1},\xi_{k+1}), \lambda_{k+1} - \lambda\ra 
&= \notag \\
&\hspace{-12em} =\alpha_{k+1}\la \nabla\Phi(\lambda_{k+1},\xi_{k+1}), \lambda_{k+1} - \zeta_k)\ra + \alpha_{k+1}\la \nabla\Phi(\lambda_{k+1},\xi_{k+1}), \zeta_k - \lambda)\ra \notag \\
 &\hspace{-12em}\stackrel{\eqref{eq:alpQuadEq},\eqref{eq:lambdakp1Def}}{=} C_{k}\la \nabla \Phi(\lambda_{k+1}, \xi_{k+1}), \eta_{k} - \lambda_{k+1}\ra + \alpha_{k+1}\la\nabla\Phi(\lambda_{k+1},\xi_{k+1}), \zeta_{k} - \lambda\ra \notag \\
  &\hspace{-12em} \stackrel{\eqref{eq_lemma1},\eqref{eq_lemma2}}{\leq }C_{k}\left(\vp(\eta_k) - \vp(\lambda_{k+1}) + \la \nabla\Phi(\lambda_{k+1},\xi_{k+1}) - \nabla \vp(\lambda_{k+1}), \eta_{k} - \lambda_{k+1}\ra \right) \notag\\
  &\hspace{-12em}+ C_{k+1}\left(  \vp(\lambda_{k+1}) - \vp(\eta_{k+1})  \right)  + V[\zeta_k](\lambda) - V[\zeta_{k+1}](\lambda) \notag\\ &\hspace{-12em}+\frac{C_{k+1}}{2L}\| \nabla\Phi(\lambda_{k+1},\xi_{k+1})- \nabla \vp(\lambda_{k+1})\|_{*}^2 \notag \\
  &\hspace{-12em}\stackrel{\eqref{eq:alpQuadEq}}{=} \alpha_{k+1} \vp(\lambda_{k+1}) + C_{k} \vp(\eta_{k}) - C_{k+1} \vp(\eta_{k+1}) + V[\zeta_k](\lambda) - V[\zeta_{k+1}](\lambda) \notag\\
   &\hspace{-12em}+C_{k}\la \nabla\Phi(\lambda_{k+1},\xi_{k+1}) - \nabla \vp(\lambda_{k+1}), \eta_{k} - \lambda_{k+1}\ra \notag \\
   &\hspace{-12em}+ \frac{C_{k+1}}{2L}\| \nabla\Phi(\lambda_{k+1}, \xi_{k+1})- \nabla\vp(\lambda_{k+1})\|_{*}^2 
\end{align*}
Rearranging terms, we obtain 
\begin{align*}
C_{k+1} \vp(\eta_{k+1}) - C_{k} \vp(\eta_{k}) 
&\leq\\
&\hspace{-11em} \leq \alpha_{k+1} \left( \vp(\lambda_{k+1}) + \la \nabla \Phi(\lambda_{k+1},\xi_{k+1}), \lambda - \lambda_{k+1}\ra \right) + V[\zeta_k](\lambda) - V[\zeta_{k+1}](\lambda)\\
&\hspace{-11em}+C_{k}\la \nabla\Phi(\lambda_{k+1}, \xi_{k+1}) - \nabla \vp(\lambda_{k+1}), \eta_{k} - \lambda_{k+1}\ra \notag \\
&\hspace{-11em}+ \frac{C_{k+1}}{2L}\| \nabla\Phi(\lambda_{k+1},\xi_{k+1})- \nabla \vp(\lambda_{k+1})\|_{*}^2. 
\end{align*}
Summing these inequalities for $k=0,\dots, N-1$, we get
\begin{align*}
C_N\vp(\eta_{N}) - C_0 \vp(\eta_{0}) 
&\leq \sum_{k=0}^{N-1} \alpha_{k+1} \left( \vp(\lambda_{k+1}) + \la \nabla\Phi(\lambda_{k+1},\xi_{k+1}), \lambda - \lambda_{k+1}\ra \right)  \\
&\hspace{-7em} + V[\zeta_0](\lambda) - V[\zeta_{N}](\lambda) + \sum_{k=0}^{N-1}C_{k}\la \nabla\Phi(\lambda_{k+1},\xi_{k+1}) - \nabla \vp(\lambda_{k+1}), \eta_{k} - \lambda_{k+1}\ra + \\
&\hspace{-7em} +\sum_{k=0}^{N-1}\frac{C_{k+1}}{2L}\|\nabla\Phi(\lambda_{k+1},\xi_{k+1})- \nabla \vp(\lambda_{k+1})\|_{*}^2.
\end{align*}
Since $C_0 = \alpha_0 = 0$ and $V[\zeta_{k}](\lambda)\geq 0$, we end up with
\begin{align*}%\label{eq:first}
C_N\vp(\eta_{N}) 
&\leq \sum_{k=0}^{N} \alpha_{k} \left( \vp(\lambda_{k}) + \la \nabla\Phi(\lambda_{k},\xi_k), \lambda - \lambda_{k}\ra \right) + V[\zeta_0](\lambda) \notag \\
&\hspace{-3em}+ \sum_{k=0}^{N-1}C_{k}\la \nabla\Phi(\lambda_{k+1}, \xi_{k+1}) - \nabla \vp(\lambda_{k+1}), \eta_{k} - \lambda_{k+1}\ra + \sum_{k=0}^{N}\frac{C_{k}}{2L}\| \nabla\Phi(\lambda_{k},\xi_k)- \nabla \vp(\lambda_{k})\|_{*}^2.
\end{align*}
Since $\lambda \in \Lambda$ was chosen arbitrarily, we take the minimum in $\lambda$ in the right hand side of this inequality and obtain the statement of the Theorem.

\subsection{Proof of Theorem \ref{Th:stoch_err}}\label{supp:the_2}
Let us introduce a set $\Lambda_R:=\{\lambda \in H^*: \|\lambda\|_2 \leq 2R\}$. Then, from \eqref{eq:ASGDConv} since $\zeta_0=0$ and  $V[\zeta](\lambda) = \frac{1}{2}\|\lambda - \zeta\|^2_2$, we have
{\small \begin{align}
C_N \vp(\eta_N) & \leq \min_{\lambda \in \dd{\Lambda}} \left\{ \sum_{k=0}^N \alpha_{k} \left( \vp(\lambda_{k}) + \la \nabla  \Phi(\lambda_{k},\xi_k), \lambda - \lambda_{k}\ra \right) + \frac12\|\lambda\|_2^2 \right\} \notag \\ 
 	&\hspace{-3em}+ \sum_{k=0}^{N-1}C_{k}\la \nabla \Phi(\lambda_{k+1},\xi_{k+1}) - \nabla \vp(\lambda_{k+1}), \eta_{k} - \lambda_{k+1}\ra + \sum_{k=0}^{N}\frac{C_{k}}{2L}\| \nabla  \Phi(\lambda_{k},\xi_k)- \nabla \vp(\lambda_{k})\|_{2}^2. \notag \\
	&\leq \min_{\lambda \in \Lambda_R} \left\{  \sum_{k=0}^N \alpha_{k} \left( \vp(\lambda_{k}) + \la \nabla  \Phi(\lambda_{k},\xi_k), \lambda - \lambda_{k}\ra \right) + \frac{1}{2} \|\lambda\|_2^2 \right\}  \notag \\ 
 	&\hspace{-3em}+ \sum_{k=0}^{N-1}C_{k}\la \nabla \Phi(\lambda_{k+1},\xi_{k+1}) - \nabla \vp(\lambda_{k+1}), \eta_{k} - \lambda_{k+1}\ra + \sum_{k=0}^{N}\frac{C_{k}}{2L}\| \nabla  \Phi(\lambda_{k},\xi_k)- \nabla \vp(\lambda_{k})\|_{2}^2 \notag \\
	&\leq \min_{\lambda \in \Lambda_R} \left\{  \sum_{k=0}^N \alpha_{k} \left( \vp(\lambda_{k}) + \la \nabla  \Phi(\lambda_{k},\xi_k), \lambda - \lambda_{k}\ra \right)  \right\} + 2R^2 \notag \\
    &\hspace{-3em}+ \sum_{k=0}^{N-1}C_{k}\la \nabla \Phi(\lambda_{k+1},\xi_{k+1}) - \nabla \vp(\lambda_{k+1}), \eta_{k} - \lambda_{k+1}\ra + \sum_{k=0}^{N}\frac{C_{k}}{2L}\| \nabla  \Phi(\lambda_{k},\xi_k)- \nabla \vp(\lambda_{k})\|_{2}^2.
\label{eq:proof_st_1}
\end{align}}
Our next goal is to take the expectation from the both sides of this inequality with respect to the seqence $\xi_0,...,\xi_N$. To do so, we iteratively, for each $j$ from $N$ to $0$ fix the history $\xi_0,...,\xi_{j-1}$ and take the expectation w.r.t $\xi_j$.

% \leq \min_{\lambda \in \Lambda} \left\{ \sum_{k=0}^N \alpha_{k} \left( \vp(\lambda_{k}) + \la \nabla  \Phi(\lambda_{k},\xi_k), \lambda - \lambda_{k}\ra \right) + V[\zeta_0](\lambda) \right\} \notag \\ 
%  &\hspace{-3em}+ \sum_{k=0}^{N-1}C_{k+1}\la \nabla \Phi(\lambda_{k+1},\xi_{k+1}) - \nabla \vp(\lambda_{k+1}), \eta_{k} - \lambda_{k+1}\ra + \sum_{k=0}^{N}\frac{C_{k}}{2L}\| \nabla  \Phi(\lambda_{k},\xi_k)- \nabla \vp(\lambda_{k})\|_{2}^2.

% From Theorem  \ref{Th:ASGDConv}, we have
% 	\begin{align}\label{eq:first}
% 	C_N\vp(\eta_N) 
%     &\leq \min_{\lambda} \left\{ \sum_{k=0}^N \alpha_{k} \left( \vp(\lambda_{k}) + \la \nabla \Phi(\lambda_{k},\xi_k), \lambda - \lambda_{k}\ra \right) + V[\zeta_0](\lambda) \right\} \notag \\ 
%  &\hspace{-4em}+ \sum_{k=0}^{N-1}C_{k+1}\la \nabla \Phi(\lambda_{k+1},\xi_k) - \nabla \vp(\lambda_{k+1}), \eta_{k} - \lambda_{k+1}\ra + \sum_{k=0}^{N}\frac{C_{k}}{2L}\| \nabla\Phi(\lambda_{k},\xi_k)- \nabla \vp(\lambda_{k})\|_{2}^2.%\label{eq:ASGDConv}
% 	\end{align}
Since $\mathbb{E}_{\xi_{k+1}}\left[ \dd{\nabla}\Phi(\lambda_{k+1},\xi_{k+1})|\xi_1,\dots,\xi_{k}\right] = \nabla \vp(\lambda_{k+1})$, $\lambda_{k+1}$ and $\eta_k$ are deterministic functions of $(\xi_1,\dots,\xi_{k})$, we have $\mathbb{E}_{\xi_1,\dots \xi_k}\la \dd{\nabla}\Phi(\lambda_{k+1},\xi_{k+1}) - \nabla \vp(\lambda_{k+1}), \eta_{k} - \lambda_{k+1}\ra = 0$. 
By the Theorem assumption
\[
\mathbb{E}_{\xi_k}\left[\left\| \nabla \Phi(\lambda_{k},\xi_k)- \nabla \vp(\lambda_{k})\right \|_2^2 |\xi_1,\dots,\xi_{k-1} \right] \leq \frac{\e L\alpha_k}{C_k}.
\]
Thus, after taking the full expectation $\E$, the last three terms in the r.h.s. of \eqref{eq:proof_st_1} satisfy
{\small \begin{align}
\label{eq:proof_st_2}
\E &\left[2R^2 + \sum_{k=0}^{N-1}C_{k+1}\la \nabla \Phi(\lambda_{k+1},\xi_{k+1}) - \nabla \vp(\lambda_{k+1}), \eta_{k} - \lambda_{k+1}\ra + \sum_{k=0}^{N}\frac{C_{k}}{2L}\| \nabla  \Phi(\lambda_{k},\xi_k)- \nabla \vp(\lambda_{k})\|_{2}^2 \right] \notag \\
&\hspace{10cm} \leq 2R^2 +\frac{C_N\e}{2},
\end{align}}
where we used that $C_N = \sum_{k=0}^N \alpha_k$.

Let us now estimate the expectation of the first term in the r.h.s. of \eqref{eq:proof_st_1}. By the definition of $F(x,\xi)$ and $F^*(-A^T\lambda,\xi)$ in subsection \ref{S:stoch_prob_setup}, we have
\begin{align}
F^*(-A^T\lambda_k,\xi_k) + \la A \nabla F^*(-A^T\lambda_k,\xi_k),\lambda_k \ra &= \la - A^T\lambda_k, x(-A^T\lambda_k,\xi_k) \ra \notag \\
& \hspace{3em} -F(x(-A^T\lambda_k,\xi_k), \xi_k) + \la Ax(-A^T\lambda_k,\xi_k), \lambda_k \ra \notag \\
&= -F(x(-A^T\lambda_k,\xi_k), \xi_k).
\end{align}
On the other hand, by Fenchel duality,
\begin{align}
\E_{\xi_k}F(x(-A^T\lambda_k,\xi_k), \xi_k) &= \E_{\xi_k} \max_{\tilde{\lambda}} \{\la x(-A^T\lambda_k,\xi_k), \tilde{ \lambda} \ra - F^*(\tilde{\lambda},\xi) \} \notag \\
& \geq \max_{\tilde{\lambda}} \{\la \E_{\xi_k}x(-A^T\lambda_k,\xi_k), \tilde{\lambda} \ra - \E_{\xi_k}F^*(\tilde{\lambda},\xi) \} = f\left(\E_{\xi_k}x(-A^T\lambda_k,\xi_k)\right).
\end{align}
Hence,
$$
\E_{\xi_k} (F^*(-A^T\lambda_k,\xi_k) + \la A \nabla F^*(-A^T\lambda_k,\xi_k),\lambda_k \ra) \leq -f(\E_{\xi_k}x(-A^T\lambda_k,\xi_k))
$$
Using this inequality, \eqref{eq:vp_def} and that $\nabla \Phi (\lambda_k,\xi_k) = b - A \nabla F^*(-A^T\lambda_k,\xi_k) = b- Ax(-A^T\lambda_k,\xi_k)$, we obtain
\begin{align}
\label{eq:proof_st_3}
\E_{\xi_k} (\vp(\lambda_{k}) + \la \nabla  \Phi(\lambda_{k},\xi_k), \lambda - \lambda_{k}\ra ) & = \la b, \lambda_k \ra + \E_{\xi_k}F^*(-A^T\lambda_k,\xi_k) + \E_{\xi_k}\la b - A \nabla F^*(-A^T\lambda_k,\xi_k), \lambda - \lambda_k \ra \notag \\
&= \E_{\xi_k}( F^*(-A^T\lambda_k,\xi_k) + \la A \nabla F^*(-A^T\lambda_k,\xi_k),\lambda_k \ra) \notag \\
& + \E_{\xi_k} \la b - Ax(-A^T\lambda_k,\xi_k), \lambda\ra \notag \\
& \leq -f(\E_{\xi_k}x(-A^T\lambda_k,\xi_k)) + \la b - A\E_{\xi_k} x(-A^T\lambda_k,\xi_k), \lambda\ra.
\end{align}
Taking the full expectation from the first term in the r.h.s. of \eqref{eq:proof_st_1} and iteratively applying \eqref{eq:proof_st_3}, we obtain
\begin{align}\label{eq:expec_min}
\E &\min_{\lambda \in \Lambda_R} \left\{  \sum_{k=0}^N \alpha_{k} \left( \vp(\lambda_{k}) + \la \nabla  \Phi(\lambda_{k},\xi_k), \lambda - \lambda_{k}\ra \right)  \right\} \leq  \min_{\lambda \in \Lambda_R} \left\{ \E \sum_{k=0}^N \alpha_{k} \left( \vp(\lambda_{k}) + \la \nabla  \Phi(\lambda_{k},\xi_k), \lambda - \lambda_{k}\ra \right)  \right\} \notag \\
& \leq \min_{\lambda \in \Lambda_R}\left\{ \sum_{k=0}^N \alpha_{k} (-f(\E x(-A^T\lambda_k,\xi_k)) +  \la b - A\E  x(-A^T\lambda_k,\xi_k), \lambda\ra ) \right\}  \notag \\
&\leq 
C_N \min_{\lambda \in \Lambda_R}\left\{ -f(\E \hat{x}_N) +  \la b - A\E \hat{x}_N, \lambda\ra \right\} \leq -C_N f(\E \hat{x}_N) + C_N \min_{\lambda \in \Lambda_R} \la b - A\E \hat{x}_N, \lambda\ra \notag \\
&= -C_N f(\E \hat{x}_N) - 2 C_N R\| b - A\E  \hat{x}_N \|_2,
\end{align}
where we also used the convexity of $f$, equality $\sum_{k=0}^N\alpha_k=C_N$, and definitions of $\hat{x}_N$ and $\Lambda_R$.

Taking the expectation in \eqref{eq:proof_st_1} and combining it with \eqref{eq:proof_st_2} and \eqref{eq:expec_min}, we obtain
\begin{align}\label{eq:func_mathExp}
	\mathbb{E}\vp(\eta_N) + f(\mathbb{E}\hat{x}_N) 
    &\leq - 2R\|A\E\hat{x}_N-b\|_2 +\frac{2R^2}{C_N} +\frac{\e}{2}.
\end{align}
Hence, by weak duality $-f(x^*)\leq \vp(\eta^*)$, 
\begin{align}\label{eq:func_math_fin}
	f(\E\hat{x}_N)  - f(x^*) \leq f(\E\hat{x}_N) + \vp(\eta^*) \leq f(\E\hat{x}_N) + \mathbb{E}\vp(\eta_N)  
    &\leq  \frac{2R^2}{C_N} +\frac{\e}{2}.
\end{align}
Since $\lm^*$ is an optimal solution of Problem (D), we have, for any $x \in Q$, $f(x^*)\leq f(x) + \la \lm^*, Ax-b \ra$.
Then using assumption $\|\lm^*\|_2\leq R$ and choosing $x = \E\hat{x}_N$, we get
\begin{equation}\label{eq:weak_dual}
f(\E\hat{x}_N) \geq f(x^*) - R\|A\E\hat{x}_N-b\|_2
\end{equation}
%This and \eqref{eq:func_mathExp}
Using this and weak duality $-f(x^*)\leq \vp(\eta^*)$ and taking the expectation, we obtain
\begin{align*}
\mathbb{E}\vp(\eta_N) +f(\mathbb{E}\hat{x}_N) \geq \vp(\eta^*) +f(\mathbb{E}\hat{x}_N) \geq -f(x^*) +f(\mathbb{E}\hat{x}_N) \stackrel{\eqref{eq:weak_dual}}{\geq} -R\|A\mathbb{E}\hat{x}_N - b\|_2
\end{align*}
Using this and \eqref{eq:func_mathExp}, we get
\begin{equation}\label{eq:MathIneq}
\|A\mathbb{E}\hat{x}_N - b\|_2 \leq \frac{2R}{C_N} + \frac{\e}{2R}
\end{equation}
It remains to estimate the growth of coefficients $C_N$. So, we prove by induction that the coefficients $C_k$ generated by Algorithm \ref{alg:main} satisfy the following condition
\begin{equation}
C_k \geq \frac{(k+1)^2}{8L}.
\label{eq:CkGrowth}
\end{equation}	
Since $C_0 = 0$ for $k=1$ $C_1  \stackrel{\eqref{eq:PDalpQuadEq}}{=} \frac{1}{2L} $ and \eqref{eq:CkGrowth} holds.
Let us now assume that \eqref{eq:CkGrowth} holds for some $k \geq 1$ and prove that it holds for $k+1$.
	By \eqref{eq:alpQuadEq}, $\alpha_{k+1}$ is the largest root of the equation
	$2L\alpha_{k+1}^2 - \alpha_{k+1} - C_k = 0.$
	Thus, 
	\begin{equation}\label{eq:alphaK}
	\alpha_{k+1} = \frac{1 + \sqrt{ 1 + 8LC_{k}}}{4L} = \frac{1}{4L} + \sqrt{\frac{1}{8L^2} + \frac{C_{k}}{2L}} \geq 
	\frac{1}{4L} + \sqrt{\frac{C_{k}}{2L}} \geq
	\frac{1}{4L} + \frac{1}{\sqrt{2L}}\frac{k+1}{2\sqrt{2L}} =
	\frac{k+2}{4L}.
	\end{equation}
	Using the induction assumption, \eqref{eq:alpQuadEq},  \eqref{eq:CkGrowth} and \eqref{eq:alphaK}, we obtain that \eqref{eq:CkGrowth} holds for $k+1$
	\begin{equation*}
	C_{k+1} = C_k + \alpha_{k+1} \geq \frac{(k+1)^2}{8L} + \frac{k+2}{4L} \geq \frac{(k+2)^2}{8L}.
	\end{equation*}
    Combining \eqref{eq:func_math_fin}, \eqref{eq:MathIneq}, and \eqref{eq:CkGrowth}, we finish our proof.

\subsection{Proof of Lemma \ref{Lm:dual_obj_properties2}}\label{supp:lemma_2}
First, let us estimate the Lipschitz constant of $\nabla \W_{\gamma}^*(\Blm)$ 
\begin{align*}
\|\nabla \W^*_\gamma(\Blm_1) - \nabla \W^*_\gamma(\Blm_2)\|^2_2
&\stackrel{\eqref{eq:DualObjGrad} }{=}\left\|
\sqrt{W} \left(
\begin{aligned}
\nabla \W^*_{\gamma,\mu_1}&([\BBlm_1]_1)\\
&...\\
\nabla \W^*_{\gamma,\mu_m}&([\BBlm_1]_m)
\end{aligned}
\right) - 
\sqrt{W} \left(
\begin{aligned}
\nabla \W^*_{\gamma,\mu_1}&([\BBlm_2]_1)\\
&...\\
\nabla \W^*_{\gamma,\mu_m}&([\BBlm_2]_m)
\end{aligned}
\right)
\right\|_2^2 \\
&\leq (\lambda_{\max}(\sqrt{W}))^2 \left\|
\begin{aligned}
\nabla \W^*_{\gamma,\mu_1}([\BBlm_1]_1) &- \nabla \W^*_{\gamma,\mu_1}([\BBlm_2]_1)\\
&...\\
\nabla \W^*_{\gamma,\mu_m}([\BBlm_1]_m) &- \nabla \W^*_{\gamma,\mu_m}([\BBlm_2]_m)
\end{aligned}
\right\|_2^2 \\
&=(\lambda_{\max}(\sqrt{W}))^2 \sum_{i=1}^m \left\| \nabla \W^*_{\gamma,\mu_i}([\BBlm_1]_i) - \nabla \W^*_{\gamma,\mu_i}([\BBlm_2]_i) \right\|_2^2\\
& \leq (\lambda_{\max}(\sqrt{W}))^2 \sum_{i=1}^m \frac{1}{\gamma^2} \left\| [\BBlm_1]_i - [\BBlm_2]_i \right\|_2^2 \\
& = \frac{(\lambda_{\max}(\sqrt{W}))^2}{\gamma^2} \sum_{i=1}^m \left\| [\sqrt{W} (\Blm_1 - \Blm_2)]_i \right\|_2^2\\
& = \frac{(\lambda_{\max}(\sqrt{W}))^2}{\gamma^2} \left\| \sqrt{W} (\Blm_1 - \Blm_2) \right\|_2^2\\
& \leq \frac{(\lambda_{\max}(\sqrt{W}))^4}{\gamma^2} \left\|\Blm_1 - \Blm_2\right\|_2^2,
\end{align*}
where we used notation $\BBlm = \sqrt{W} \Blm$, the definition of matrix $\sqrt{W}$, $1/\gamma$-Lipschitz continuity of $\nabla \W_{\gamma,\mu_i}^*(\blm_i)$ for all $i=1,\dots, m$. Since $(\lambda_{\max}(\sqrt{W}))^4 = (\lambda_{\max}(W))^2$, we obtain that the dual function $ \W_{\gamma}^*(\Blm)$ has $\lambda_{\max}(W) /\gamma$-Lipschitz continuous gradient.  

By Lemma \ref{Lm:dual_obj_properties}, vectors $p_j(\blm_j,Y_r^j)$, $j=1,...,m$, $r=1,...,M$ defined in \eqref{eq:StochGradDef} satisfy $\E_{Y_r^j} p_j(\blm_j,Y_r^j) = \nabla\W^*_{\gamma,\mu_j}(\blm_j)$. Thus, by \eqref{eq:DualObjGrad}, \eqref{eq:tnW_def}, \eqref{eq:StochGradDef0} we have $\E \widetilde{\nabla}\W_{\gamma}^*(\Blm) = \nabla\W_{\gamma}^*(\Blm)$.

Further, for $j=1,...,m$, we estimate the variance of $p_j(\blm_j, Y^j)$
 \begin{align*} 
\mathbb{E}_{Y^j\sim \mu_j} \| p_j(\blm_j, Y^j) - \nabla \mathcal{W}^*_{\gamma,\mu_j}(\blm_j)\|^2_2&=
\mathbb{E}_{Y^j\sim \mu_j} \sum_{l=1}^n \left(\frac{\exp\left(([\blm_j]_l -c_l(Y))/\gamma\right)}{\sum^n_{\ell=1}\exp\left(([\blm_j]_\ell -c_\ell(Y))/\gamma\right)} - [\nabla \mathcal{W}^*_{\gamma,\mu_j}(\blm_j)]_l\right)^2 \\
&=\sum_{l=1}^n \mathbb{E}_{Y^j\sim \mu_j}\frac{\exp^2\left(([\blm_j]_l -c_l(Y))/\gamma\right)}{\left(\sum^n_{\ell=1}\exp\left(([\blm_j]_\ell -c_\ell(Y))/\gamma\right)\right)^2} - \sum_{l=1}^n  [\nabla \mathcal{W}^*_{\gamma,\mu_j}(\blm_j)]^2_l \\
&\leq \sum_{l=1}^n\int_{\mathcal{Y}} \frac{\exp^2(([\blm_j]_l-c_l(y))/\gamma) q(y) }{\left(\sum_{\ell=1}^n\exp(([[\blm]_j]_\ell-c_\ell(y))/\gamma)\right)^2} q(y) dy \\
&=\int_{\mathcal{Y}} \frac{\sum_{l=1}^n\exp^2(([\blm_j]_l-c(y,z_i))/\gamma)  }{\left(\sum_{\ell=1}^n\exp(([[\blm]_j]_\ell-c_\ell(y))/\gamma)\right)^2} q(y)dy \leq \int_{\mathcal{Y}} q(y)dy = 1.
\end{align*}
Hence, by \eqref{eq:StochGradDef0}, for $j=1,...,m$, we have
\begin{equation}\label{eq:var-M}
\E_{Y_r^j\sim\mu_j, r=1,...,M} \|\widetilde{\nabla} \W_{\gamma,\mu_j}^*(\blm_j) - \nabla \W_{\gamma,\mu_j}^*(\blm_j)\|_2^2 \leq \frac{1}{M}.
\end{equation}

% Furthermore, we estimate the parameter of Lipschitz continuity of $\nabla \W_{\gamma}^*(\Blm)$   
% %$\nabla \W^*_{\gamma}(\sqrt{W}y) = \sum_{i,j=1}^m\sqrt{W}_{ij}\W^*_{\gamma,\mu_j}(\lambda_j)$ as 
% %as $L = \|\sqrt{W}\|^2_{2 \to 2}/\gamma = \lambda_{\max}(W) /\gamma$~\cite{Kakade2009}, where $\lambda_{\max}(W)$ is the maximum eigenvalue of the matrix $W$.
% \begin{align*}
% \|\nabla \W^*_\gamma(\sqrt[]{W}\Blm) - \nabla \W^*_\gamma(\sqrt[]{W}\Bmu)\|^2_2
% &=\|\sqrt{W} \nabla \W^*_\gamma(\BBlm) - \sqrt{W}\nabla \W^*_\gamma(\BBmu) \|_2^2 \leq \|\sqrt{W} \|^2_2\|\nabla \W^*_\gamma(\BBlm) -\nabla \W^*_\gamma(\BBmu) \|_2^2 \\
% &= \lm^2_{\max}(W) \sum_{i=1}^m\|\nabla \W^*_\gamma(\blm) -\nabla \W^*_\gamma(\bmu)\|^2_2 \leq \frac{\lm^2_{\max}(W) }{\gamma^2} \sum_{i=1}^m \|\blm_i - \bmu_i\|^2_2 \\
% &= \frac{\lm^2_{\max}(W) }{\gamma^2} \|\BBlm - \BBmu\|^2_2 = \frac{\lm^2_{\max}(W) }{\gamma^2} \|\sqrt[]{W}\Blm - \sqrt[]{W}\Bmu\|^2_2 
% \end{align*}
% where in last inequality we used $1/\gamma$-Lipschitz continuity of $\nabla \W_{\gamma}^*(\blm_i)$ for all $i=1,\dots, m$. Hence, dual function $ \W_{\gamma}^*(\Blm)$ has $\lambda_{\max}(W) /\gamma$-Lipschitz continuous gradient.   

By the same arguments as above for the estimate of the Lipschitz constant for $\nabla \W_{\gamma}^*(\Blm)$, we estimate the variance of $\widetilde{\nabla} \W_{\gamma}^*(\Blm)$. Denoting $\E = \E_{Y_r^j\sim\mu_j, j=1,...,m, r=1,...,M}$, we have 
\begin{align*}
\E
\|\widetilde{\nabla} \W_{\gamma}^*(\Blm) - \nabla \W_{\gamma}^*(\Blm)\|_2^2 &\stackrel{\eqref{eq:DualObjGrad},\eqref{eq:tnW_def}}{=} \E \left\|
\sqrt{W} \left(
\begin{aligned}
\widetilde{\nabla} \W^*_{\gamma,\mu_1}&([\BBlm]_1)\\
&...\\
\widetilde{\nabla} \W^*_{\gamma,\mu_m}&([\BBlm]_m)
\end{aligned}
\right) - 
\sqrt{W} \left(
\begin{aligned}
\nabla \W^*_{\gamma,\mu_1}&([\BBlm]_1)\\
&...\\
\nabla \W^*_{\gamma,\mu_m}&([\BBlm]_m)
\end{aligned}
\right)
\right\|_2^2\\
& \leq (\lambda_{\max}(\sqrt{W}))^2 \E\left\|
\begin{aligned}
\widetilde{\nabla} \W^*_{\gamma,\mu_1}([\BBlm]_1) &- \nabla \W^*_{\gamma,\mu_1}([\BBlm]_1)\\
&...\\
\widetilde{\nabla} \W^*_{\gamma,\mu_m}([\BBlm]_m) &- \nabla \W^*_{\gamma,\mu_m}([\BBlm]_m)
\end{aligned}
\right\|_2^2
\\
& = (\lambda_{\max}(\sqrt{W}))^2 \E\sum_{i=1}^m \left\| \widetilde{\nabla}\W^*_{\gamma,\mu_i}([\BBlm]_i) - \nabla \W^*_{\gamma,\mu_i}([\BBlm]_i) \right\|_2^2
\\
& \stackrel{\eqref{eq:var-M}}{\leq} \frac{(\lambda_{\max}(\sqrt{W}))^2m}{M} = \frac{\lambda_{\max}(W)m}{M},
\end{align*}
% \begin{align*}
% &E_{Y_r^j\sim\mu_j, j=1,...,m, r=1,...,M} 
% \|\widetilde{\nabla} \W_{\gamma}^*(\sqrt[]{W}\Blm) - \nabla \W_{\gamma}^*(\sqrt[]{W}\Blm)\|_2^2 \\
% &\hspace{3em}= \|\sqrt[]{W}\|^2_2 E_{Y_r^j\sim\mu_j, j=1,...,m, r=1,...,M} \|\widetilde{\nabla} \W_{\gamma}^*(\BBlm) - \nabla \W_{\gamma}^*(\BBlm)\|_2^2\\
% &\hspace{3em}= \lm^2_{\max}(W)E_{Y_r^j\sim\mu_j, j=1,...,m, r=1,...,M} \sum_{i=1}^m \|\widetilde{\nabla} \W_{\gamma, \mu_i}^*(\blm_i) - \nabla \W_{\gamma,\mu_i}^*(\blm_i)\|_2^2\\
% &\hspace{3em}\leq \lm^2_{\max}(W) \sum_{j=1}^m E_{Y_r^j\sim\mu_j, r=1,...,M}\|\widetilde{\nabla} \W_{\gamma, \mu_j}^*(\blm_j) - \nabla \W_{\gamma,\mu_j}^*(\blm_j)\|_2^2 
% \stackrel{\eqref{eq:var-M}}{\leq} \frac{\lm^2_{\max}(W)m}{M},
% \end{align*}
% \begin{align*}
% \E_{Y_r^j\sim\mu_j, j=1,...,m, r=1,...,M} \|\widetilde{\nabla} \W_{\gamma}^*(\sqrt[]{W}\Blm) - \nabla \W_{\gamma}^*(\sqrt[]{W}\Blm)\|_2^2 \leq \frac{\lambda_{\max}(W)}{M}, \; \Blm \in \R^{mn},
% \end{align*}
which finishes the proof of the Lemma.

\subsection{Proof of Theorem \ref{Th:WBCompl}}\label{supp:the_3}
Combining Lemma \ref{Lm:dual_obj_properties2} and Theorem \ref{Th:stoch_err} 
for our particular case of primal-dual pair of problems \eqref{consensus_problem2}-\eqref{eq:DualPr}
with $A=\sqrt{W}$, $b=0$, $L= \lambda_{\max}(W)/\gamma$, since $N=\sqrt{32\lambda_{\max}(W)R^2/(\e\gamma)}$, we obtain the first statement of the theorem.

Let us now estimate the overall complexity of Algorithm \ref{alg:main}. For each agent $i$, the complexity of each iteration is dominated by the complexity of calculation of stochastic approximation $\widetilde{\nabla} \W_{\gamma,\mu_i}^*([\BBlm_{k}]_i)$ for the gradient. This complexity is $O(mnM_k)$. Thus, to get the overall complexity, we need to estimate $\sum_{k=1}^{N}M_k$ 
\begin{align*}
\sum_{k=1}^{N}M_k &= \sum_{k=1}^{N}\max\left\{ 1,~ \left\lceil \frac{m\gamma C_{k}}{\alpha_k\e} \right\rceil \right\} \stackrel{\eqref{eq:PDalpQuadEq}}{\leq} \max\left\{N,~  \left\lceil \frac{2\lambda_{\max}(W)m}{\e}\sum_{k=1}^{N}\alpha_k\right\rceil \right\} \notag \\
& = \max\left\{N,  ~\left\lceil\frac{2\lambda_{\max}(W)m}{\e}C_N\right\rceil \right\}
 %\leq \max\left\{\sqrt{\frac{\eig R^2}{\e \gamma}}, ~ \frac{\lambda_{\max}(W)m R^2}{\e^2}\right\}, 
\end{align*}
where we used that $\sum_{k=1}^N\alpha_k=C_N$.
\dd{
From \eqref{eq:MathIneq} and definition of $N$ it follows that
\begin{align*}
\frac{2R}{C_N} \leq \frac{\e}{2R} \quad \mbox{and} \quad \frac{2R}{C_{N-1}} \geq \frac{\e}{2R}.
\end{align*}
Then 
\begin{align}\label{eq:constC}
C_{N-1} \leq 4R^2/\e 
\end{align}
From \eqref{eq:alphaK}
\begin{align}\label{eq:alpha_1}
\alpha_N = \frac{1}{4L} + \sqrt{\frac{1}{8L^2} + \frac{C_{N-1}}{2L}} \leq \frac{1}{2L} + \sqrt{\frac{C_{N-1}}{2L}} \stackrel{\eqref{eq:PDalpQuadEq}}{=}  \frac{1}{2L}+ \alpha_{N-1}
\end{align}
On the other hand, from \eqref{eq:alphaK} it follows that
\begin{align}\label{eq:alpha_2}
\alpha_{N-1} = \frac{1}{4L} + \sqrt{\frac{1}{8L^2} + \frac{C_{N-2}}{2L}} \geq \frac{1}{4L} 
\end{align}
Hence, from \eqref{eq:alpha_1} and \eqref{eq:alpha_2} we have
\begin{align*}
\alpha_N \leq 2\alpha_{N-1} + \alpha_{N-1} = 3\alpha_{N-1}\stackrel{\eqref{eq:PDalpQuadEq}}{\leq} 3C_{N-1}
\end{align*}
Since this inequality and \eqref{eq:PDalpQuadEq} we obtain $C_N \leq 4 C_{N-1}$. Then using \eqref{eq:constC} we have 
\begin{align}\label{ineq:C_N}
\sum_{k=1}^{N}M_k \leq \max\left\{\sqrt{\frac{32\eig R^2}{\e \gamma}}, ~ \frac{32\lambda_{\max}(W)m R^2}{\e^2}\right\},
\end{align}
}
where in last equality we used $N = \sqrt[]{32\lambda_{\max}(W)R^2/(\e\gamma)}$. 
To obtain the total complexity, we multiply the above estimate for $\sum_{k=1}^{N}M_k$ by $mn$. 
% For the second variants of $M_k$ we have 
% \[Q =\sum_{i=0}^{N}M_k =N = n^2\sqrt{16L R^2/\e } = \sqrt{16\eig R^2/(\e \gamma)}.\]
% Considering the complexity of calculating the gradient $\nabla \W_{\gamma}^*(\Blm)$, which is $O(n)$, we can estimate the overall complexity of the Algorithm as 
% \begin{align*}
% S \approx \max\left\{n ~\ceil{8\lambda_{\max}(W)R^2/\e}, n\sqrt{16\eig R^2/(\e \gamma)} \right\}.
% \end{align*}

\newpage

\subsection{Visualization of the Network Topologies used in Simulations}\label{supp:networks}

\begin{figure}[H]
	\centering
    \subfigure[Star Graph]{\includegraphics[origin=c,width=0.3\textwidth]{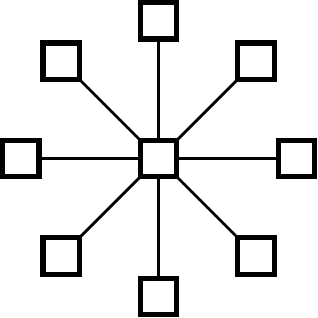}}
	\subfigure[Cycle Graph]{\includegraphics[origin=c,width=0.3\textwidth]{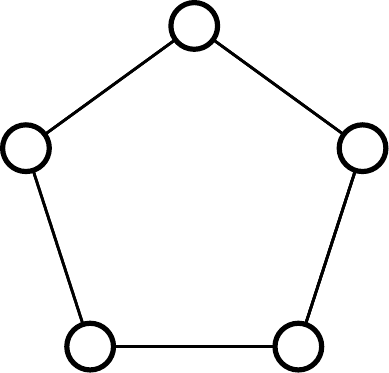}}
	\subfigure[Erd\H{o}s-R\'enyi random graph]{\includegraphics[origin=c,width=0.3\textwidth]{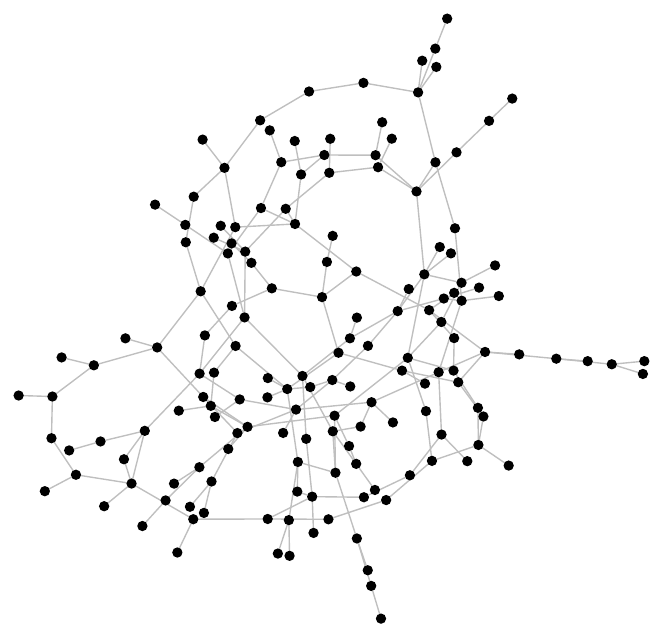}}
    \subfigure[Complete Graph]{\includegraphics[origin=c,width=0.3\textwidth]{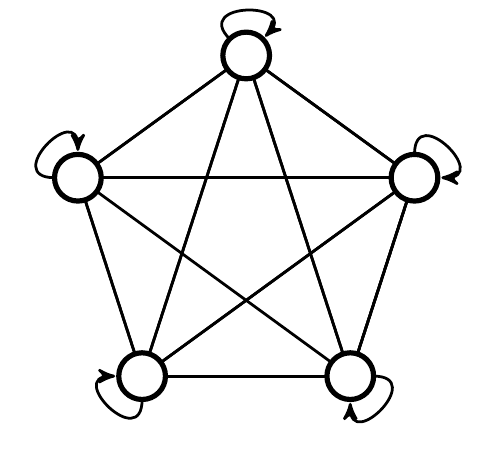}}
\caption{ Example of Network Topologies.}
    \label{fig:networks}
\end{figure}

\newpage

\subsection{Additional Example on Gaussian Distributions}\label{supp:gauss}

\begin{figure}[H]
	\centering
    \includegraphics[width=1\textwidth]{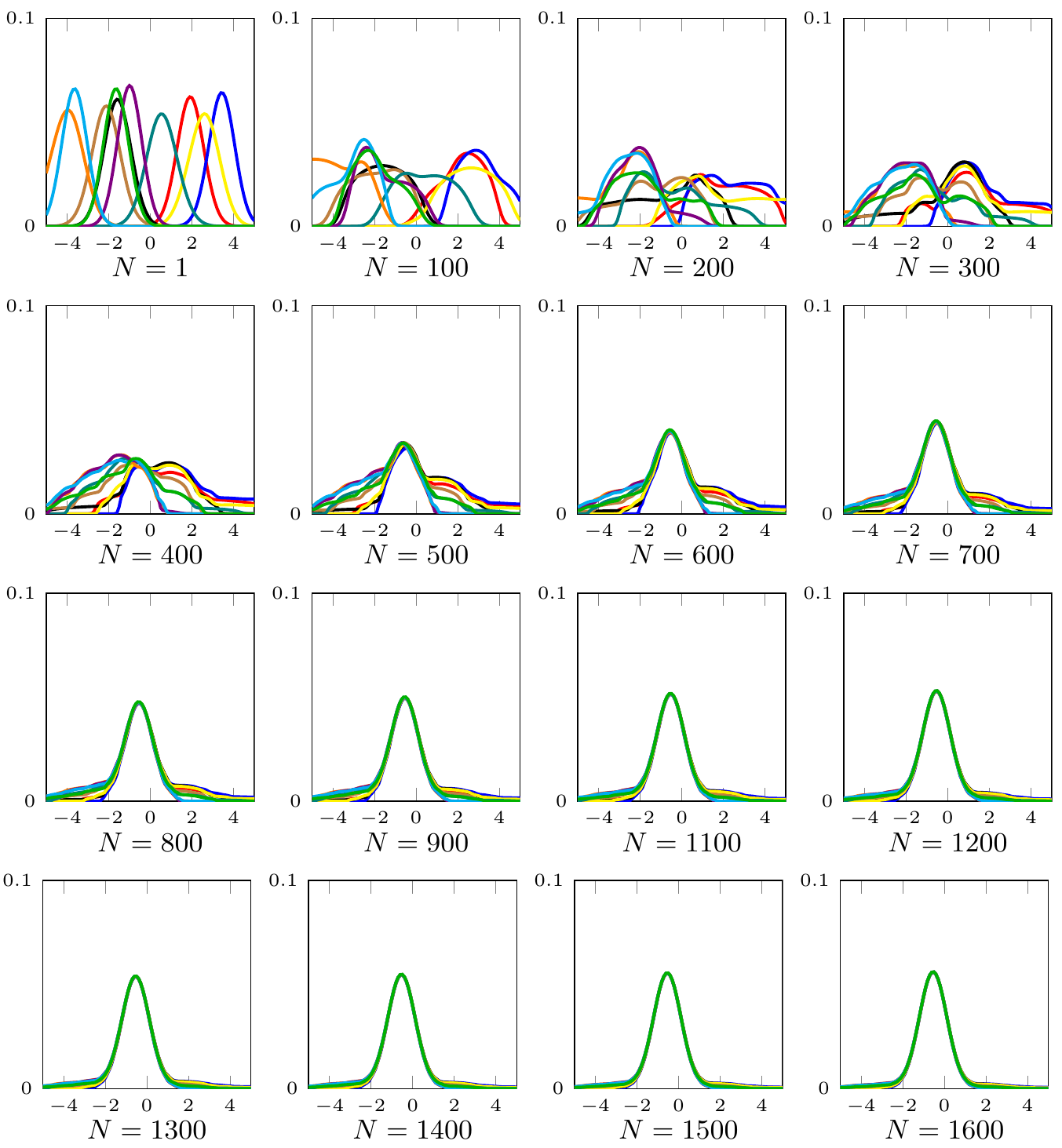}
\caption{Local Wasserstein barycenter of $10$ agents connected on an Erd\H{o}s-R\'enyi random graph. Each agent holds a private Gaussian measure from which it can query samples. Different colors represent different agents.}
    \label{fig:gauss_complete}
\end{figure}

\newpage 
\subsection{Additional Example on von Mises Distributions}\label{supp:von}

\begin{figure}[H]
	\centering
    \includegraphics[width=1\textwidth]{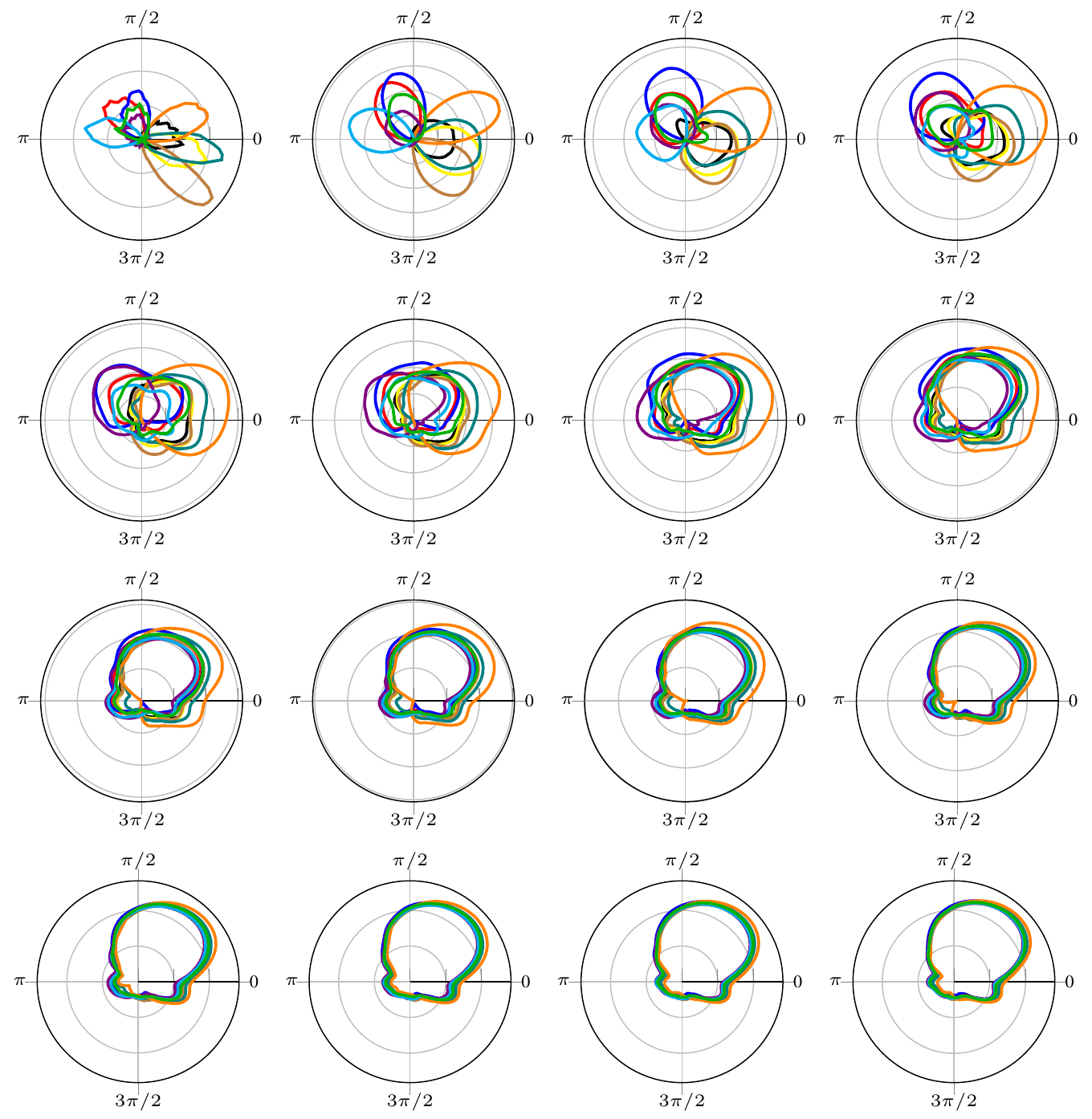}
\caption{Local Wasserstein barycenter of $10$ agents connected on an Erd\H{o}s-R\'enyi random graph. Each agent holds a private von Mises measure from which it can query samples. Different colors represent different agents. }
    \label{fig:von_complete}
\end{figure}

\newpage

\subsection{Additional Information for the MNIST Dataset}\label{supp:mnist}

\begin{figure}[H]
	\centering
    {\includegraphics[origin=c,width=0.5\textwidth]{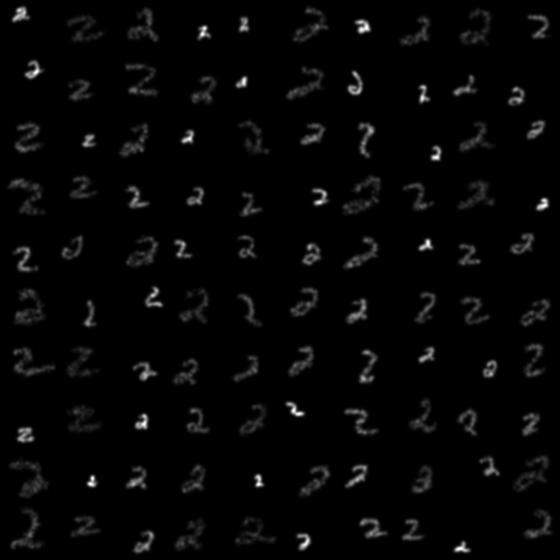}}\\
	{\includegraphics[origin=c,width=0.3\textwidth]{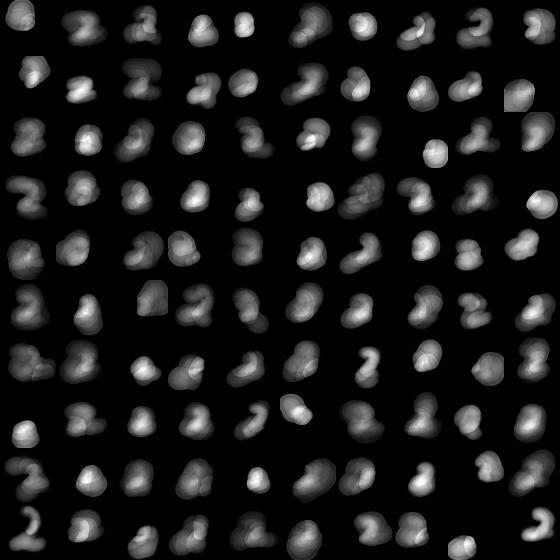}}
	{\includegraphics[origin=c,width=0.3\textwidth]{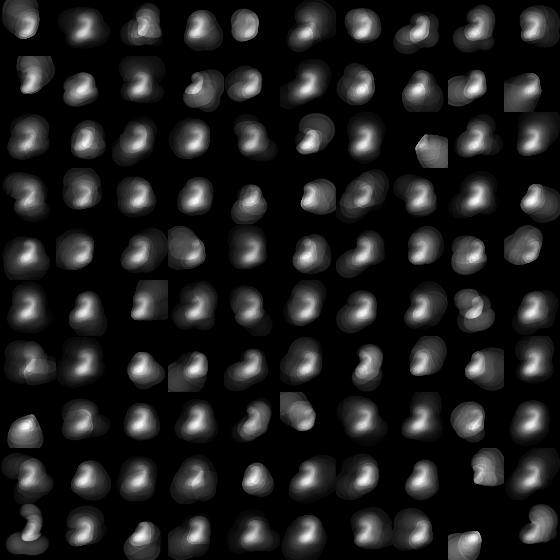}}
	{\includegraphics[origin=c,width=0.3\textwidth]{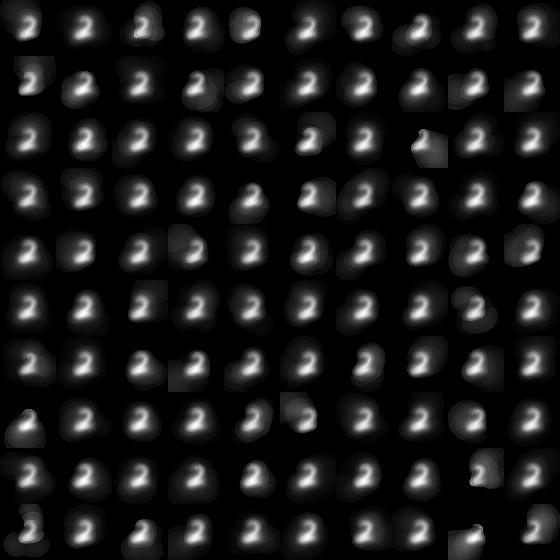}}
    {\includegraphics[origin=c,width=0.3\textwidth]{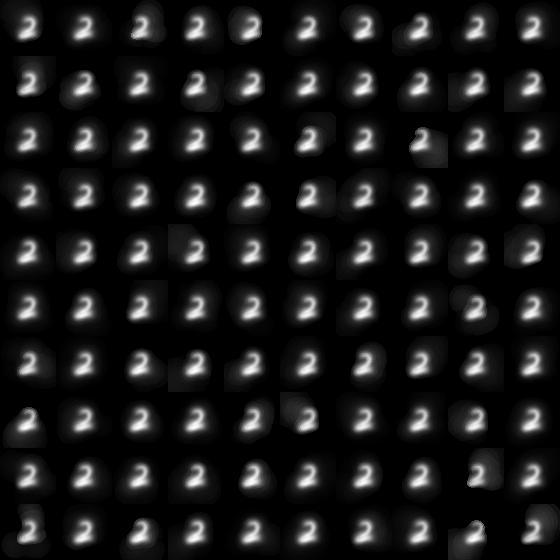}}
    {\includegraphics[origin=c,width=0.3\textwidth]{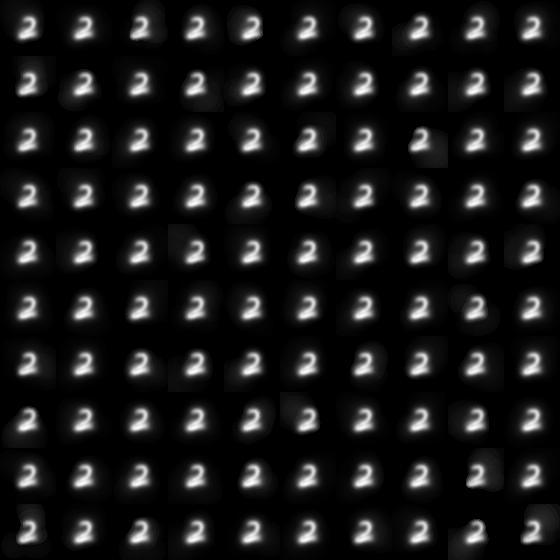}}
    {\includegraphics[origin=c,width=0.3\textwidth]{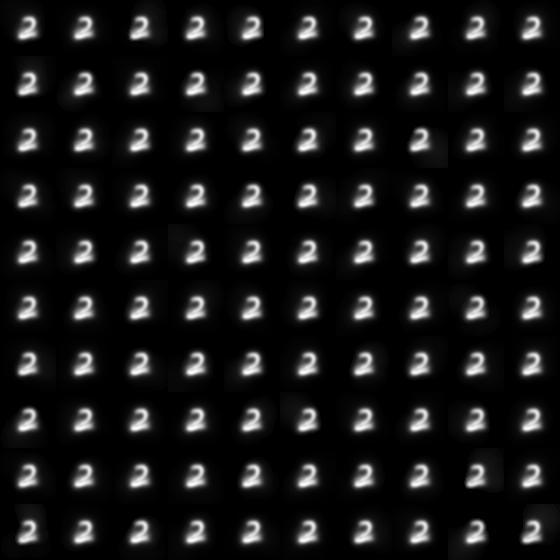}}\\
\caption{Local Wasserstein barycenter of $100$ agents connected on an Erd\H{o}s-R\'enyi random graph. Each agent holds a private sample of the digit $2$ from the MNIST dataset. We assume the normalize image as a probability distribution from which agents can sample from. }
    \label{fig:mnist_complete}
\end{figure}

\newpage

\subsection{Additional Information for the IXI Dataset}\label{supp:ixi}

\begin{figure}[H]
	\centering
    {\includegraphics[angle=180,origin=c,width=0.23\textwidth]{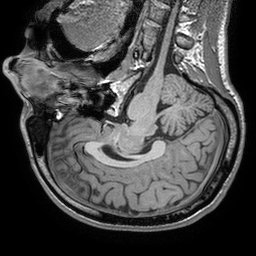}}
	{\includegraphics[angle=180,origin=c,width=0.23\textwidth]{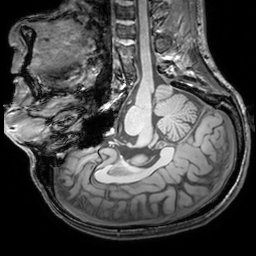}}
	{\includegraphics[angle=180,origin=c,width=0.23\textwidth]{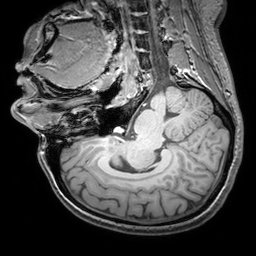}}
	{\includegraphics[angle=180,origin=c,width=0.23\textwidth]{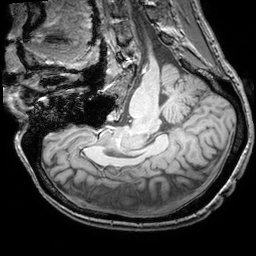}}
\caption{The samples from the IXI dataset held by four agents.}
    \label{fig:MRI_original}
\end{figure}

\begin{figure}[H]
	\centering
    {\includegraphics[origin=c,width=0.3\textwidth]{im_1e4_med2_1}}
	{\includegraphics[origin=c,width=0.3\textwidth]{im_1e4_med2_101}}
	{\includegraphics[origin=c,width=0.3\textwidth]{im_1e4_med2_1001}}
	{\includegraphics[origin=c,width=0.3\textwidth]{im_1e4_med2_8711}}
    {\includegraphics[origin=c,width=0.3\textwidth]{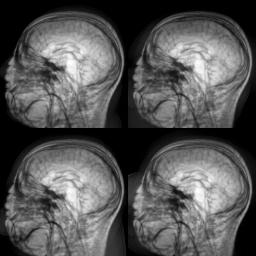}}
    {\includegraphics[origin=c,width=0.3\textwidth]{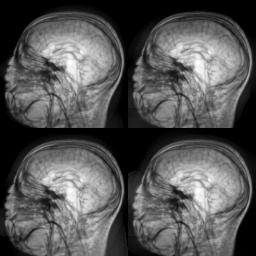}}
\caption{Local Wasserstein barycenter of $4$ agents connected on a cycle graph. Each agent holds a private sample of an magnetic resonance image from the IXI dataset. We assume the normalize image as a probability distribution from which agents can sample from. Time evolves with the number of iterations.}
    \label{fig:MRI_complete}
\end{figure}

\end{document}